\documentclass{article}
\usepackage{array}
\usepackage{amsmath}
\usepackage{amssymb}
\usepackage{bbm}
\usepackage{amsfonts}
\usepackage{graphicx}
\usepackage{caption}
\usepackage{url}
\usepackage{maplestd2e}

\DefineParaStyle{Maple Heading 1}
\DefineParaStyle{Maple Text Output}
\DefineParaStyle{Maple Dash Item}
\DefineParaStyle{Maple Bullet Item}
\DefineParaStyle{Maple Normal}
\DefineParaStyle{Maple Heading 4}
\DefineParaStyle{Maple Heading 3}
\DefineParaStyle{Maple Heading 2}
\DefineParaStyle{Maple Warning}
\DefineParaStyle{Maple Title}
\DefineParaStyle{Maple Error}
\DefineCharStyle{Maple Hyperlink}
\DefineCharStyle{Maple 2D Math}
\DefineCharStyle{Maple Maple Input}
\DefineCharStyle{Maple 2D Output}
\DefineCharStyle{Maple 2D Input}

\usepackage[top=.9in,bottom=.9in]{geometry}
 \numberwithin{equation}{section}
 
\DefineParaStyle{Maple Output}
\DefineCharStyle{2D Comment}
\DefineCharStyle{2D Math}
\DefineCharStyle{2D Output}
\begin{document}
\pagestyle{plain}
\newcommand{\ignore}[1]{}
\numberwithin{equation}{section}
\begin{flushleft} \vskip 0.3 in  
\centerline{\bf Uncovering functional relationships at zeros with special reference to Riemann's Zeta Function} \vskip .3in
\centerline{ M.L. Glasser\footnote{laryg@clarkson.edu}}\vskip .2in

\centerline{Department of Physics, Clarkson University}
 
\centerline{Potsdam, NY 13699-5820}
\centerline{}
\centerline{and}\vskip .2in
\centerline{ Michael Milgram\footnote{mike@geometrics-unlimited.com}}

\centerline{Consulting Physicist, Geometrics Unlimited, Ltd.}
\centerline{Box 1484, Deep River, Ont. Canada. K0J 1P0}
\centerline{}
\centerline{May 18, 2015}
\centerline{}
\vskip .5in
\centerline{\bf Abstract}\vskip .1in

A Master equation has been previously obtained which allows the analytic integration of a fairly large family of functions provided that they possess simple properties. Here, the properties of this Master equation are explored, by extending its applicability to a general range of an independent parameter. Examples are given for various values of the parameter using Riemann's Zeta function as a template to demonstrate the utility of the equation. The template is then extended to the derivation of various sum rules among the zeros of the Zeta function as an example of how similar rules can be obtained for other functions.
\vskip .1in

\section{Introduction}
In a previous work \cite{Master}, several ``Master'' equations were presented that permit the analytic integration of suitably restricted, but otherwise arbitrary, functions. In general, any integral of the form
\begin{equation}
\int_{-\infty}^{\infty}\,f(v)\,dv
\end{equation}

where $f(v)$ is suitably bounded as $v\rightarrow \pm\infty$ and satisfies
\begin{equation}
f(v)+f(-v-i)=0
\label{f(v)}
\end{equation}
can be written in terms of the sum of its residues residing in the region 
\begin{equation}
\nonumber
\mathfrak{S} \equiv -1 < \mathfrak{Im}(v) < 0 \,.
\end{equation} 

For the remainder of this paper, we define
\begin{equation}
{
z=\mapleinline{inert}{2d}{4*a*v*(v+I)}
{$\displaystyle 4\,av \left( v+i \right) $}\,;
}
\label{z_def}
\end{equation}
and, in addition, use $h,j,k,m,n,p,q \in \mathbf{N}$ unless otherwise specified, and set $0\leq w\leq1\,,\, w \in \mathfrak{R} $. Assuming that $F(z)$ is regular in $\mathfrak{S}$, for small, real values of the parameter $a>0$, a very simple example of such a Master equation is

\begin{equation}
\mapleinline{inert}{2d}{Int(w^(v*(v+I))*F(4*a*v*(v+I))/cosh(Pi*v), v = -infinity .. infinity) = w^{1/4}*F(a)}
{\[\displaystyle \int _{-\infty }^{\infty }\!{\frac {{w}^{v \left( v+i \right) } 
F \left( z  \right) }{\cosh \left( \pi \, 
\mbox{}v \right) }}{dv}=w^{1/4}F \left( a \right) \]}.
\label{Master1}
\end{equation}

The main motivation for this work is to investigate  the conditions under which poles of the function $F(z)$ in the complex $z$-plane will contribute to poles in $\mathfrak{S}$ in the complex $v$-plane by choosing as an appropriate example, the function $F(z)=1/\zeta(z)$, Riemann's zeta function (e.g. \cite{Edwards}). The goal is to identify terms that must be included in the right-hand side of \eqref{Master1} for various choices of the parameter $a$, and, at the same time, demonstrate how to uncover new functional relationships that may exist between the integral and residue representations inherent in \eqref{Master1} for any particular choice of $F(z)$. 
\newline

In \cite{Master}, the variable $a$ was treated as a real, positive parameter; here we intend to investigate the effect of extending $a$ to the complex domain and, from that basis, following \cite{GlasserSumRule} we further show how to develop functional relationships among complex zeros of any function F(z), using the Riemann function $\zeta(z)$ as an interesting, and challenging, example.
\newline

Fundamental to the analysis is the requirement that properties of $F(z)$ first be known in $\mathfrak{S}$, which in turn suggests that the geometry of the mapping \eqref{z_def} be examined. This is done in Section 2 for various general ranges of real and imaginary values of the parameter $a$. Armed with an understanding of the nature of \eqref{z_def}, it is then possible to judiciously translate contours of integration, making allowance for residues associated with poles for a particular choice of $F(z)$, wherever they may lie. This is done in Section 3 for special choices of the parameter $a$. In Section 4, we broaden our horizons by choosing more varied examples for $F(z)$, restricted to variants of $\zeta(z)$ in order to demonstrate the type of results that can emerge for other choices of $F(z)$. The choice $F(z)=1/\zeta(z)$ was motivated by general interest in the properties of the complex zeros of $\zeta(z)$, since the complexity and distribution of it's zeros allows us to demonstrate many of the intricacies that can be dealt with using the methods established here. The results obtained here differ from those usually found in the literature because they generally  demonstrate a relationship among function values at the zeros of $\zeta(z)$, rather than a relationship that may exist among the values of the zeros themselves as is usual (e.g. \cite[Eqs. 3.2(7) and 3.8(4)]{Edwards}, \cite{Cerone},\cite{Guillera},\cite{Matiyasevich}). This comment likely applies to any other choice of $F(z)$ \cite{deLyra},\cite{Baricz}.

\subsection{Notation and Assumptions}

Since we will be utilizing the properties of Riemann's $\zeta$ function intensively, for the sake of completeness, we remind the reader that $\zeta(z)=0$ when $z=-2n$ (so-called trivial zeros) and $z=\tau$ (non-trivial zeros). The general non-trivial zeros $\zeta(\tau)=0$ are indexed by $\tau=\sigma+i\rho$ without reference to location within  the critical strip, defined by $0<\sigma<1$ with $\rho\geq0$. When referring to zeros on the critical line, we write $\tau=\frac{1}{2}+i\rho$ in general ($\rho>0$), and  $\tau_n=\frac{1}{2}+i\rho_n$ in particular. When used in this form to index a sum over all non-trivial zeros, the implication is that $\tau$ spans all possible zeros with which it is associated, those being $\tau, \overline{\tau}, 1-\tau$ and $1-\overline{\tau}$. In some cases, the indexing over non-trivial zeros is constrained to a particular range, and this will be denoted by constraints on $\tau$ (e.g. $\Im(\tau)>0$). Throughout, we also follow the usual assumption that zeros of $\zeta(\tau)$ are simple (e.g. \cite{Bui}); consequently we always employ $\zeta^{\prime}(\tau)\neq 0$.
\newline

Finally, the appearance of $\zeta^{\prime}(-2\,n)$ is frequently replaced by a well-known identity \cite[Eq. 25.6.13]{NIST}  
\begin{equation}
\mapleinline{inert}{2d}{Zeta(1,-2*n) = 1/2*(-1)^n*Zeta(2*n+1)*(2*n)!/((2*Pi)^(2*n));}{%
\[
\zeta^{\prime} ( - 2\,n)={\displaystyle \frac {1}{2}} \,
{\displaystyle \frac {(-1)^{n}\,\zeta (2\,n + 1)\,(2\,n)\mathrm{!
}}{(2\,\pi )^{2\,n}}\,.} 
\]
}
\label{ZetaprimeNegN}
\end{equation}

\section{The nature of the transform \eqref{z_def}}

\begin{figure}
\includegraphics[width=0.5\textwidth]{{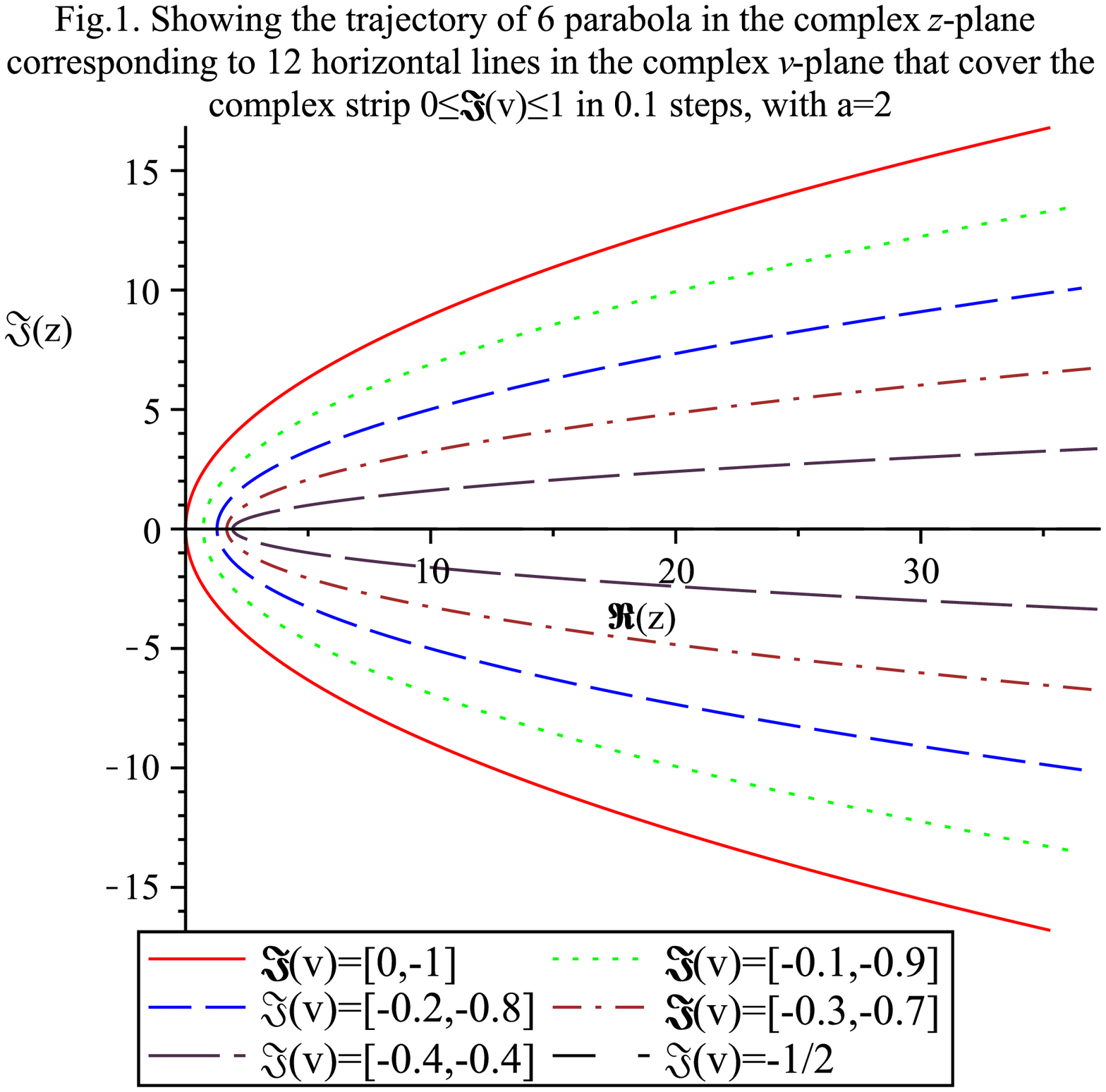}}
\caption{\label{fig:Fig1}}
\end{figure}

Figure 1 illustrates how horizontal lines spanning $\mathfrak{S}$ in the $v$-plane (separated by units of $0.1$ along the $\Im(v)$ axis) transform into parabolas in the complex $z$-plane if $a \in \Re, a>0$. The range $\Re(v)<0$ corresponds to the lower branch of each of the curves, the range $\Re(v)>0$ to the upper branch. In other words, the integration region $\mathfrak{S}$ upon which \eqref{Master1} is based transforms into the interior of the parabolic region bounded by the solid curve in Figure \ref{fig:Fig1}. If F(z) is singularity free in this region, \eqref{Master1} is valid. Otherwise modifications must be made. In the remainder of this work, we will be considering the case where F(z) is replaced by its reciprocal, so \eqref{Master1} becomes
\begin{equation}
\mapleinline{inert}{2d}{Int(w^(v*(v+I))*F(4*a*v*(v+I))/cosh(Pi*v), v = -infinity .. infinity) = w^{1/4}*F(a)}
{\[\displaystyle \int _{-\infty }^{\infty }\!{\frac {{w}^{v \left( v+i \right) } 
}{F \left( z \right) \cosh \left( \pi \, 
\mbox{}v \right) }}{dv}=\frac{w^{1/4}}{F \left( a \right)} \]}
\label{1/Master1}
\end{equation}
subject to the requirement that F(z) has no zeros in $\mathfrak{S}$. This simple strategy allows us to relate residues of $F(z)$ if zeros of $F(z)$ do happen to arise in $\mathfrak{S}$, to values of the same function as determined by the parameter $a$. As stated, in most of the remainder of this work we shall, in \eqref{1/Master1}, use $F(z)=\zeta(z)$. 
\newline

With reference to Figure \ref{fig:Fig1} where we have used $a=2$, as the magnitude of $a$ increases, ($a \in\mathfrak{R},a>0$) the width of the parabolic region increases, and eventually swallows the first known complex conjugate zero pair of $\zeta(z)$ at $z=1/2 \pm i\rho_{1}$ ($\rho_{1}=14.134725142...$) when $a = \rho_{1}^{2}/2\approx 100$. In the complex $v$-plane, for the case $z=1/2 + i\rho_{1}$ corresponding poles appear at complex points $v = (\frac{1}{2\rho_{1}},0)$ and
$v=(-\frac{1}{2\rho_{1}},-1)$. For the associated complex conjugate zero at $z=1/2 - i\rho_{1}$, the same occurs when $v = (-\frac{1}{2\rho_{1}},0)$ and $v=(\frac{1}{2\rho_{1}},-1)$. Thus two complex poles, corresponding to two complex zeros of $\zeta(z)$ in the $z$ plane, are reflected in \eqref{1/Master1} by four complex poles of $1/\zeta(z)$ in the $v$ plane  and so \eqref{1/Master1} must be modified accordingly when $a \geq \rho_{1}^{2}/2$. This will be explored in the next section.
\newline

In the case that $a \in \Re, a<0$, the parabola in Figure \ref{fig:Fig1} becomes its mirror image reflected about the vertical axis, and so the parabolic region in the $z$-plane corresponding to $\mathfrak{S}$ in the $v$-plane always encloses the `trivial' zeros of $\zeta(z)$ at $z=-2n$. In the case that $n < -a/2$, the pole belonging to $1/\zeta(z)$ at $z=-2n$ generates two sets of poles in the $v$-plane at 
\begin{equation}
\mapleinline{inert}{2d}{V[p1] := 0, -1/2+(1/2)*sqrt(1+2*n/a)}
{\[v_{-}= \,(0,\,-\frac{1}{2}\pm \frac{1}{2}\, \sqrt{1+\,{\frac {2n}{a}}})\]}
\end{equation}
and for the case $n>-a/2 $, two similar sets of poles appear at 
\begin{equation} 
\mapleinline{inert}{2d}{V[m1] := -(1/2)*sqrt(-1-2*n/a), -1/2}
{\[v_{+}=( \,\pm\frac{1}{2}\, \sqrt{-1-\,{\frac {2n}{a}}},\,-\frac{1}{2})\]}.
\end{equation}
To generalize, when $a<0$, the result \eqref{1/Master1} must always be corrected for the presence of singularities, if $F(z)$ has zeros anywhere on the negative real $z$-axis.
\newline

In the case that $a =i\beta, \beta\in \Re$, the family of parabolas in Figure \ref{fig:Fig1} rotates, opening upward ($\beta>0$) or downward ( $\beta<0$). In these two cases, corrections must be made for complex poles lying in the upper (lower) half $z$-plane respectively, but not in the form of complex conjugate pairs. Specifically, for the case $\beta>0$, as the magnitude of $\beta$ increases from zero, the parabola opens wider to envelop more (lower-lying) poles; the final singularity in the $z$-plane occurs when the upward facing parabola envelopes the first (lowest) zero of $\zeta(z)$ at $z=1/2 + i\rho_{1}$. This happens when $\beta$ increases through the value $1/(16\rho_{1})$, yielding the last entry of a family of poles in the $v$-plane terminating at the points
\begin{equation}
v_{1}=(2\,\rho_{1},-1)
\end{equation}
and
\begin{equation}
v_{2}=(-2\,\rho_{1},0).
\end{equation}

A similar result holds for the case $\beta<0$ where the conjugate zero $z=1/2 - i\rho_{1}$  generates two similar points in the $v$-plane. Notice that in these cases, conjugate zeros of $\zeta(z)$ can never contribute in pairs, but individually, each pole belonging to a zero in the $z$-plane that does contribute, will generate two contributing poles in the $v$-plane. In the general case, these observations would apply to any function $F(z)$ that might possess conjugate zeros in the upper and lower $z$ half-plane.
\newline

Finally, we turn our attention to the general case $a =\alpha+i\beta$, and consider a variety of cases that depend on the magnitude of $\alpha$ relative to $\beta$ relative to the location of the pole of $F(z)$ in the $z$-plane. With reference to Figure \ref{fig:Fig1}, the parabolic region of interest will be oriented off the coordinate axes, depending on the sign of $\alpha$ and $\beta$ and their relative magnitudes, but the vertex remains at the origin.
Consider \eqref{1/Master1}, and suppose
\begin{equation}
F(\sigma + i \rho) =0\;.
\end{equation}
There are several cases to consider.

\subsection{Case 1: $\alpha\rho+(\sigma-2\alpha)\beta>0$}

If the parabolic boundary of the region in the $z$-plane defined by some choice of $\alpha+i\beta$ passes through the complex point $(\sigma,\rho)$, then two corresponding points exist in the complex $v$-plane. The first, defined by $(v_{1,R},v_{1,I})$  is given by co-ordinates 
\begin{equation}
\mapleinline{inert}{2d}{v[R] = (1/4)*(a*sqrt(2*X+2*Y)+b*sqrt(2*X-2*Y))/(a^2+b^2)}{\[\displaystyle v_{{1,R}}={\frac {\alpha \sqrt{2\,X+2\,Y}+\beta \sqrt{2\,X-2\,Y}}{4({\alpha}^{2}+{\beta}^{2})}}\]}
\end{equation}
\begin{equation}
\mapleinline{inert}{2d}{V[I] = -1/2(1/4)*(a*sqrt(2*X-2*Y)-b*sqrt(2*X+2*Y))/(a^2+b^2)}{\[\displaystyle v_{{1,I}}=-\frac{1}{2}+{\frac {\alpha \sqrt{2\,X-2\,Y}-\beta \sqrt{2\,X+2\,Y}}{4({\alpha}^{2}+{\beta}^{2})}}\]}
\end{equation}
where
\begin{equation}
\mapleinline{inert}{2d}{X = sqrt((a^2+b^2)*((a-sigma)^2+(b-rho)^2))}{\[\displaystyle X= \sqrt{ \left( {\alpha}^{2}+{\beta}^{2} \right)  \left(  \left( \alpha-\sigma \right) ^{2}+ \left( \beta-\rho \right) ^{2} \right) }\]}
\end{equation}
and
\begin{equation}
\mapleinline{inert}{2d}{Y = b^2-a^2-b*rho+a*sigma}{\[\displaystyle Y={\beta}^{2}-{\alpha}^{2}-\beta\rho+\alpha\sigma\]}.
\end{equation}
It is worth noting that, in the real plane, these co-ordinates satisfy 
\begin{equation}
\mapleinline{inert}{2d}{v[R]^2+(v[I]+1/2)^2 = (1/4)*sqrt(((\alpha-sigma)^2+(\beta-rho)^2)/(\alpha^2+\beta^2))}{\[\displaystyle {v_{{1,R}}}^{2}+ \left( v_{{1,I}}+1/2 \right) ^{2}=\frac{1}{4} \sqrt{{\frac { \left( \alpha-\sigma \right) ^{2}+ \left( \beta-\rho \right) ^{2}}{{\alpha}^{2}+{\beta}^{2}}}}\]}
\label{VRI_sum}
\end{equation}
and
\begin{equation}
\mapleinline{inert}{2d}{v[R]^2-(v[I]+1/2)^2 = -1/4+(1/4)*\alpha*sigma/(\alpha^2+b^2)+(1/4)*rho*b/(\alpha^2+b^2)}{\[\displaystyle {v_{{1,R}}}^{2}- \left( v_{{1,I}}+1/2 \right) ^{2}=-\frac{1}{4}+{\frac {\alpha\sigma}{4({\alpha}^{2}+{\beta}^{2})}}+{\frac {\beta\rho}{4({\alpha}^{2}+{\beta}^{2})}}\]}\,,
\label{VRI_diff}
\end{equation}
leading to
\begin{equation}
\mapleinline{inert}{2d}{(v[I]+1/2)^2 = (1/8)*sqrt(((\alpha-sigma)^2+(b-rho)^2)/(\alpha^2+b^2))+(1/8)*(\alpha*(\alpha-sigma)+b*(b-rho))/(\alpha^2+b^2)}{\[\displaystyle  \left( v_{{1,I}}+1/2 \right) ^{2}=\frac{1}{8} \sqrt{{\frac { \left( \alpha-\sigma \right) ^{2}+ \left( \beta-\rho \right) ^{2}}{{\alpha}^{2}+{\beta}^{2}}}}+{\frac {\alpha \left( \alpha-\sigma \right) +\beta \left( \beta-\rho \right) }{8\,({\alpha}^{2}+{\beta}^{2})}}\]}\,.
\end{equation}
The condition that the point $(v_{1,R},v_{1,I})$ in the complex $v$-plane lies inside $\mathfrak{S}$ is that the imaginary coordinate $ v_{1,I}$ satisfies $-1\leq v_{1,I} \leq 0$, equivalent to the requirement that
\begin{equation}
\mapleinline{inert}{2d}{(v[I]+1/2)^2 = (1/8)*sqrt(((a-sigma)^2+(b-rho)^2)/(a^2+b^2))+(1/8)*(a*(a-sigma)+b*(b-rho))/(a^2+b^2)}{\[\displaystyle  \left( -1/2 \right) ^{2}\leq \frac{1}{8} \sqrt{{\frac { \left( \alpha-\sigma \right) ^{2}+ \left( \beta-\rho \right) ^{2}}{{\alpha}^{2}+{\beta}^{2}}}}+{\frac {\alpha \left( \alpha-\sigma \right) +\beta \left( \beta-\rho \right) }{8\,({\alpha}^{2}+{\beta}^{2})}}\]}\leq ({1}/{2})^{2}\,.
\label{CondCase1}
\end{equation}

As noted, this case possesses a second solution defined by the complex point $(v_{2,R},v_{2,I})$ given by
\begin{eqnarray}
v_{2,R}=&-v_{1,R} \nonumber\\
v_{2,I}=&-v_{1,I}-{\displaystyle 1}\,
\label{v2}
\end{eqnarray}
which shares the same condition \eqref{CondCase1} in order that this point lies inside $\mathfrak{S}$.

\subsection{Case 2: $\alpha\rho+(\sigma-2\alpha)\beta<0$}

In this case, the corresponding complex points  $({v}_{1,R},{v}_{1,I})$ in the $v$-plane are given by
\begin{equation}
\mapleinline{inert}{2d}{UC1R := (1/4)*(a*sqrt(2*X+2*Y)-b*sqrt(2*X-2*Y))/(a^2+b^2)}
{\[\displaystyle {{v}_{1,R}}= \,{\frac {\alpha \sqrt{2\,X+2\,Y}-\beta \sqrt{2\,X-2\,Y}}{4({\alpha}^{2}+{\beta}^{2})}}\]}
\end{equation}
and
\begin{equation}
\mapleinline{inert}{2d}{UC1I := -1/2-(1/4)*(\alpha*sqrt(2*X-2*Y)+\beta*sqrt(2*X+2*Y))/(a^2+b^2)}{\[\displaystyle  {{v}_{1,I}}= \,-\frac{1}{2}-{\frac {\alpha \sqrt{2\,X-2\,Y}+\beta \sqrt{2\,X+2\,Y}}{4({\alpha}^{2}+{\beta}^{2})}}\]},
\end{equation}
both of which obey the equivalent of \eqref{VRI_sum} and \eqref{VRI_diff}, and hence the condition \eqref{CondCase1} holds in this case. The second point belonging to this case obeys \eqref{v2}

\subsection{Case 3: $\alpha\rho+(\sigma-2\alpha)\beta=0$}

Solving for $\beta$ we have
\begin{equation}
\beta=\alpha\rho/(2\alpha-\sigma)
\end{equation}
giving, for the case $\alpha<0, \alpha-\sigma<0$
\begin{equation}
v_{1,R}=\frac{\rho\sqrt{\sigma-\alpha}}{2\sqrt{-\alpha}\sqrt{(2a-\sigma)^2+\rho^2}}
\end{equation}
\begin{equation}
\mapleinline{inert}{2d}{-1/2+(1/2)*sqrt(-\alpha)*(-2*\alpha+sigma)*sqrt(sigma-\alpha)/(\alpha*sqrt(4*\alpha^2-4*\alpha*sigma+sigma^2+rho^2))}{\[\displaystyle v_{1,I}= -1/2+1/2\,{\frac { \sqrt{-\alpha} \left( -2\,\alpha+\sigma \right)  \sqrt{\sigma-\alpha}}{\alpha \sqrt{{(2\alpha-\sigma)}^{2}+{\rho}^{2}}}}\]}
\end{equation}
or, if $\alpha>0, \alpha-\sigma>0$
\begin{equation}
v_{1,R}=\frac{\rho\sqrt{\alpha-\sigma}}{2\sqrt{\alpha}\sqrt{(2a-\sigma)^2+\rho^2}}
\end{equation}
\begin{equation}
\mapleinline{inert}{2d}{-1/2+(1/2)*sqrt(-\alpha)*(-2*\alpha+sigma)*sqrt(sigma-\alpha)/(\alpha*sqrt(4*\alpha^2-4*\alpha*sigma+sigma^2+rho^2))}{\[\displaystyle v_{1,I}= -1/2+1/2\,{\frac {  \left( 2\,\alpha-\sigma \right)  \sqrt{\alpha-\sigma}}{\sqrt{\alpha} \sqrt{({2\alpha-\sigma})^{2}+{\rho}^{2}}}}\]}
\end{equation}
If $\alpha>0, \sigma>\alpha $ we find
\begin{equation}
\mapleinline{active}{2d}{-(2*a-sigma)*sqrt(-((2*a-sigma)^2+rho^2)*(a-sigma))/(2*sqrt(a)*((2*a-sigma)^2+rho^2))}{\[v_{1,R}=\pm\,{\displaystyle {\frac { \left( 2\,a-\sigma \right)  \sqrt{- \left(  \left( 2\,\alpha-\sigma \right) ^{2}+{\rho}^{2} \right)  \left(\alpha -\sigma \right) }}{2 \sqrt{\alpha} \left(  \left( 2\,\alpha-\sigma \right) ^{2}+{\rho}^{2} \right)} 
\mbox{}}}\]}
\end{equation}
\begin{equation}
\mapleinline{active}{2d}{-1/2+(1/2)*rho*sqrt(-((2*\alpha-sigma)^2+rho^2)*(-sigma+\alpha))/(sqrt(\alpha)*((2*\alpha-sigma)^2+rho^2))}{\[v_{1,I}=-1/2\mp{\displaystyle \frac {\rho\, \sqrt{- \left(  \left( 2\,\alpha-\sigma \right) ^{2}+{\rho}^{2} \right)  \left( \alpha-\sigma \right) }}{2 \sqrt{\alpha} \left(  \left( 2\,\alpha-\sigma \right) ^{2}+{\rho}^{2} \right) 
\mbox{}}}\]}
\end{equation}
according as $\sigma<2a$ or $\sigma>2a$ respectively. For completeness sake, we include the case $\alpha<0, \sigma<\alpha$, which does not apply to the choice $F(z)=\zeta(z)$ since it is known that $\zeta(\sigma+i\rho)$ does not possess complex zeros when $\sigma<0$. The corresponding results are

\begin{equation}
\mapleinline{active}{2d}{-(2*\alpha-sigma)*sqrt(-((2*\alpha-sigma)^2+rho^2)*(\alpha-sigma))/(2*sqrt(\alpha)*((2*\alpha-sigma)^2+rho^2))}{\[v_{1,R}=\mp\,{\displaystyle {\frac { \left( 2\,\alpha-\sigma \right)  \sqrt{ \left(  \left( 2\,\alpha-\sigma \right) ^{2}+{\rho}^{2} \right)  \left(\alpha -\sigma \right) }}{2 \sqrt{-\alpha} \left(  \left( 2\,\alpha-\sigma \right) ^{2}+{\rho}^{2} \right)} 
\mbox{}}}\]}
\end{equation}
\begin{equation}
\mapleinline{active}{2d}{-1/2+(1/2)*rho*sqrt(-((2*\alpha-sigma)^2+rho^2)*(-sigma+\alpha))/(sqrt(\alpha)*((2*\alpha-sigma)^2+rho^2))}{\[v_{1,I}=-1/2\mp{\displaystyle \frac {\rho\,\sqrt{-\alpha} \sqrt{ \left(  \left( 2\,\alpha-\sigma \right) ^{2}+{\rho}^{2} \right)  \left( \alpha-\sigma \right) }}{2 \alpha \left(  \left( 2\,\alpha-\sigma \right) ^{2}+{\rho}^{2} \right) 
\mbox{}}}\]}
\end{equation}
again, according as $\sigma<2a$ or $\alpha>\sigma>2a$ respectively. In all these cases, a second set of solutions exists. Each set satisfies \eqref{v2}.

\section{Consequences}
\subsection{The case $0<a<\rho_{1}^{2}/2, a \in \Re$}
With the caveats given and, because we have chosen $F(z)=\zeta(z)$ and $\zeta(z)$ contains no zeros in $\mathfrak{S}$ for $0<a<\rho_{1}^{2}/2$ (see Section 2),  \eqref{1/Master1} becomes

\begin{maplegroup}
\mapleresult
\begin{maplelatex}
\begin{equation}
\mapleinline{inert}{2d}{Int(w^(v*(v+I))*Pi/(Zeta(4*a*v*(v+I))*cosh(Pi*v)), v = -infinity .. infinity) = Pi*w^(1/4)/Zeta(a)}{\[\displaystyle \int _{-\infty }^{\infty }\!{\frac {{w}^{v \left( v+i \right) } }{\zeta  \left( 4\,av \left( v+i \right)  \right) 
\mbox{}\cosh \left( \pi \,v \right) }}{dv}={\frac { \,w^{1/4}}{\zeta  \left( a \right) }}\,\hspace{2cm 0\leq w\leq1\,}.\]} 
\label{Master1_Zeta}
\end{equation}
\end{maplelatex}
\end{maplegroup}

A further variation can be obtained by shifting the contour of integration in \eqref{Master1_Zeta} upwards (parallel to the real $v$ axis). In terms of Figure (\ref{fig:Fig1}) the associated boundary contour in the $z$-plane moves to the left and opens to swallow poles belonging to both trivial zeros on the real axis, and non-trivial zeros in the critical strip as $\Im(v)$ increases. 
\newline

The evaluation of \eqref{Master1_Zeta} during translation, is achieved by identifying residues of known singularities at $v=(k+1/2)\,i,\,\,$ $z=-2n$ and $z=\tau$, encountered as the contour progresses through the upper half of the $v$-plane. Noting that the integral itself vanishes on the boundary $\Im(v)\rightarrow \infty$, we find, for $0\leq w\leq1\,$, the sum rule \cite{GlasserSumRule} 

%
%

\begin{equation}
\mapleinline{inert}{2d}{Sum(w^(tau/a)/(-a+4*tau)^(1/2)/Zeta(1,tau)/sinh(1/2*Pi*(-1+4*tau/a^(1
/2))^(1/2)),tau) =
-1/2*a^(1/2)*w^(1/4)/Pi/Zeta(1/4*a)+Sum(w^(-2*n/a)/(a+8*n)^(1/2)/Zeta(
1,-2*n)/sin(1/2*Pi*(1+8*n/a^(1/2))^(1/2)),n = 1 ..
infinity)+a^(1/2)*Sum((-1)^k*w^(-1/4*(2*k+3)*(2*k+1))/Zeta(-1/4*a*(2*k
+3)*(2*k+1)),k = 0 .. infinity)/Pi;}{%
\maplemultiline{
{\displaystyle \sum _{\tau \geq \tau_{_{ M}} }} \,{\displaystyle \frac {w^{ {
\tau }/{4\,a}}}{\sqrt{ - a + \tau\, }\,\zeta^{\prime} (\tau )\,\mathrm{
sinh} \Big(   {\displaystyle \frac {\pi \,\sqrt{ - 1 + 
{\displaystyle  {\tau/ }{{a}}} }}{2}} \Big) }
} =  -{\displaystyle \frac {2\,\sqrt{
a}\,w^{1/4}}{\pi \,\zeta ({\displaystyle  {a}} )}} 
\mbox{}  +\,S_{1} \,+\, S_{2} \\ }
}
\label{Srule1}
\end{equation}
\newline

where
\begin{equation}
S_{1} \equiv  {\displaystyle \sum _{n=1}^{\infty }} \,
{\displaystyle \frac {w^{ - {n}/{2\,a}}}{\sqrt{a + 2\,n}\,
\zeta^{\prime} ( - 2\,n)\,\mathrm{sin} \Big( {\displaystyle 
\frac {\pi \,\sqrt{1 + {\displaystyle {2\,n}/{{a}}} }}{
2}} \, \Big) }}  
\mbox{}
\label{S1}
\end{equation}
\begin{equation}
S_{2}\equiv -{\displaystyle \frac {4\,\sqrt{a}}{\pi }\,  
{\displaystyle \sum _{k=1}^{\infty }} \,{\displaystyle \frac {(-1
)^{k}\,w^{ 1/{4}-k^2}}{\zeta (  
{\displaystyle {a\,(1-4\,k^{2})}} )}}  }
\label{S2}
\end{equation}

and we reiterate that ${\displaystyle \sum_{\tau\geq \tau_{_M}}}$ represents the sum over all specified values of $\tau_{m}$ for which $\zeta(\tau_{m})=0$ with $m\geq M$ and for the moment we take $M=1$ (see subsection \ref{Case4}). Each of the three summations in \eqref{Srule1} corresponds to the sum over the named residues attached to each set of poles with indices labelled corresponding to the discussion above.

\subsection{Case $a>\rho^{2}_{M}/2$, $a \in \Re$} \label{Case4}

As indicated in Section 2, if $\rho^{2}_{2}>2\,a>\rho^{2}_{1}$ the bounding parabola in the $z$-plane encloses the pole of $\zeta(\rho_{1})$, thereby generating a pole in the region $\mathfrak{S}$ and violating the premise of \eqref{1/Master1}. To compensate, the residue of each such pole must be added to the right-hand side of \eqref{Master1_Zeta}; similarly, if $\rho^{2}_{M+1}>2\,a>\rho^{2}_{M}$, (where $M$ labels a non-trivial zero of $\zeta$ according to the order in which they will encounter the parabola in the $z$-plane as $a$ increases in magnitude), then the sum of the residues of all $M$ poles encountered must be added to the right-hand side of \eqref{Master1_Zeta} in the form of a finite sum of M terms equivalent to the left-hand side of \eqref{Srule1}. When the integration contour is moved upwards according to the derivation of \eqref{Srule1} new residues are encountered, such that a sum over all non-trivial zeros in \eqref{Srule1} effectively only includes contributions from zeros lying outside the parabola. These contributions in turn, are dictated by the magnitude of the parameter $a$; this is indicated by setting $M>1$ appropriately in \eqref{Srule1}.        

\subsection{Exceptional cases with $a>0$ \label{Sub_exc}}

Examination of \eqref{S1} indicates that scattered terms of  the sum diverge when 
\begin{equation}
1+2\,n/a=4\,h^2\,\,;
\label{C1}
\end{equation}
similarly, scattered terms of \eqref{S2} diverge when
\begin{equation}
a(4k^2-1)=2j
\label{C2}
\end{equation}
where $h$ and $j$ are positive integers. Both these Diophantine conditions are simultaneously fulfilled when $a>0$ is of the form
\begin{equation}
{\displaystyle a_{p,q}=\frac{2\,p}{4q^2-1}}\,.
\label{apq}
\end{equation}
where $p$ and $q$ are positive integers. Specific diverging terms occur when the indices $n$ and $k$ obey
\begin{eqnarray} \label{n&k}
n=&j   \nonumber \\
k=&h-1. 
\end{eqnarray}

Taking the limit $a\rightarrow a_{p,q}$ shows that the $n^{th}$ indexed exceptional term of \eqref{S1} diverges as
\begin{eqnarray}
\mapleinline{inert}{2d}{1/(w^(1/2*n/a)*(a+2*n)^(1/2)*Zeta(1,-2*n)*sin(1/2*Pi*(1+2*n/a)^(1/2))
) =
-2*(-1)^h*a^(3/2)*w^(-1/2*n/a)/n/Zeta(1,-2*n)/Pi/(a-apq)-(-1)^h*w^(-1/
2*n/a)*ln(w)/a^(1/2)/Pi/Zeta(1,-2*n)-1/4*(-1)^h*(a+3*n)*w^(-1/2*n/a)/a
^(1/2)/n/h^2/Pi/Zeta(1,-2*n);}{%
\maplemultiline{
{\displaystyle \frac {w^{(-\frac {n}{2\,a})}}{\,\sqrt{a + 2\,n}\,
\zeta^{\prime}( - 2\,n)\,\mathrm{sin} \Big(  \! {\displaystyle 
\frac {\pi \,\sqrt{1 + {\displaystyle  {2\,n}/{a}} }}{2}} 
 \!  \Big) }} \approx \\
 - {\displaystyle \frac {2\,(-1)^{h}\,a^{(3/2)}\,w^{( - \frac {n
}{2\,a})}}{n\,\zeta^{\prime}( - 2\,n)\,\pi \,(a - \mathit{a_{p,q}})}} 
 - {\displaystyle \frac {(-1)^{h}\,w^{( - \frac {n}{2\,a})}\,
\mathrm{ln}(w)}{\sqrt{a}\,\pi \,\zeta^{\prime}( - 2\,n)}}  - 
{\displaystyle \frac {1}{4}} \,{\displaystyle \frac {(-1)^{h}\,(a
 + 3\,n)\,w^{( - \frac {n}{2\,a})}}{\sqrt{a}\,n\,h^{2}\,\pi \,
\zeta^{\prime}( - 2\,n)}}  }
}
\label{Except1}
\end{eqnarray}

and the $k^{th}$ indexed exceptional term of \eqref{S2}, diverges as

\begin{equation}
\mapleinline{inert}{2d}{4*a^(1/2)*(-1)^k/Pi/(w^(1/4*(2*k+3)*(2*k+1)))/Zeta(-a*(2*k+3)*(2*k+1)
) =
2*(-1)^(k+1)*a^(3/2)*w^(-1/2*j/a)/j/Pi/Zeta(1,-2*j)/(a-apq)-a^(1/2)*(-
1)^k*w^(-1/2*j/a)/j/Pi/Zeta(1,-2*j)-2*Zeta(2,-2*j)*a^(1/2)*(-1)^k*w^(-
1/2*j/a)/Pi/Zeta(1,-2*j)^2;}{%
\maplemultiline{
{\displaystyle \frac {4\,\sqrt{a}\,(-1)^{k}}{\pi \,w^{(\frac {(2
\,k + 3)\,(2\,k + 1)}{4})}\,\zeta ( - a\,(2\,k + 3)\,(2\,k + 1))}
} \approx \\
{-\displaystyle \frac {2\,(-1)^{(k)}\,a^{(3/2)}\,w^{( - 
\frac {j}{2\,a})}}{j\,\pi \,\zeta^{\prime} (- 2\,j)\,(a - \mathit{
a_{p,q}})}}  - {\displaystyle \frac {\sqrt{a}\,(-1)^{k}\,w^{( - 
\frac {j}{2\,a})}}{j\,\pi \,\zeta^{\prime} (- 2\,j)}}  - 
{\displaystyle \frac {2\,\zeta^{\prime\prime} (- 2\,j)\,\sqrt{a}\,(-1)^{k}
\,w^{( - \frac {j}{2\,a})}}{\pi \,\zeta^{\prime} (- 2\,j)^{2}}}  }
}
\label{Except2}
\end{equation}
\newline

Following \eqref{n&k}, we identify the index $n$ of \eqref{Except1} with the parameter $j$ of \eqref{Except2}, and the index $k$ of \eqref{Except2} with the parameter $h-1$ of \eqref{Except1}, and discover that divergent terms of two different sums cancel pairwise and corresponding terms of the two sums reduce to the $0^{th}$ order (non-divergent) term(s) given in the respective equations \eqref{Except1} and \eqref{Except2} whenever the corresponding summation index in either sum satisfies \eqref{C1} or \eqref{C2}. The general form of \eqref{Srule1} is fairly lengthy in this case, but can be written with recourse to \eqref{Except1} and \eqref{Except2}. This is left as an exercise for the reader.
\newline

As an example, in the case $a\rightarrow a_{p,q}=2/35$ ($p=1,q=3$), divergent terms in the sum \eqref{S1} are indexed by $n=1,33,37,117...$ and the corresponding cancelling divergent terms of the sum \eqref{S2} are indexed by $k=2,16,17,31...$. In terms of the geometry of Figure (\ref{fig:Fig1}), this case arises for special values of $a$ (see \eqref{apq}) when, first of all, poles of $1/\zeta(z)$ belonging to $z=-2n$ and poles of $1/\cosh(\pi v)$ at $v=(2k+1)\,i\,/2$ simultaneously lie on the integration contour in the $v$-plane. The residue singularity becomes a dipole as the integration contour is shifted upward (see the derivation of \eqref{Srule1}) whenever these two poles occasionally coalesce, and the corresponding scattered divergences in \eqref{Srule1} reflect this eventuality. 
\newline

A further solution of \eqref{C1} and \eqref{C2} exists when $a=2m, m=1,2,...$. A careful derivation of this case following the above analysis by evaluating the limits as $a\rightarrow2m$ yields the results \eqref{Except1} and \eqref{Except2} with $a \rightarrow 2m$ except that divergent terms in both results are each multiplied by a factor of two; they still cancel pairwise.
For  the sake of brevity, we give here the sum rule corresponding to \eqref{Srule1} when $a=2m$ and $w=1$:
\begin{equation}
\mapleinline{inert}{2d}{Sum(1/((-2*m+tau)^(1/2)*Zeta(1,tau)*sinh(1/4*Pi*(-4+2*tau/m)^(1/2))),
tau) =
-2*m^(1/2)/Pi/Zeta(2*m)+Sum(1/((m+n)^(1/2)*Zeta(1,-2*n)*sin(1/2*Pi*(1+
n/m)^(1/2))),n =
1)-1/8*2^(1/2)*m^(1/2)*Pi*Sum((-1)^h*(-1+12*h^2)/n/h^2/Zeta(1,-2*n),h)
-2^(1/2)*m^(1/2)/Pi*Sum((-1)^k/Zeta(1,-2*j)*(1/j+2*Zeta(2,-2*j)/Zeta(1
,-2*j)),k = 0 .. infinity);}{%
\maplemultiline{
{\displaystyle \sum _{\tau }} \,{\displaystyle \frac {1}{\sqrt{
 - 2\,m + \tau }\,\zeta^{\prime} (\tau )\,\mathrm{sinh} \big( 
{\displaystyle \frac {\pi}{2}} \,\sqrt{ - 1 + {\displaystyle{\tau }/{2m}} }  \;  \big) }} =   
\mbox{} {\displaystyle \frac{1}{\sqrt{\,2}}} {\displaystyle \sideset{}{\,^{'}}\sum_{n=1}} \,
{\displaystyle \frac {1}{\sqrt{m + n}\,\zeta^{\prime}( - 2\,n)\,
\mathrm{sin} \Big(  \! {\displaystyle \frac {\,\pi \,\sqrt{1 + 
{\displaystyle  {n}/{m}} }}{2}}   \Big) }}\\   
\mbox{} - {\displaystyle \frac {\pi}{8}} \,\sqrt{2\,m} 
\,  {\displaystyle \sum _{n=1}}^{\star} \,{\displaystyle \frac {(
-1)^{h}\,( - 1 + 12\,h^{2})}{n\,h^{2}\,\zeta^{\prime} ( - 2\,n)}} 
\mbox{} - {\displaystyle  { \frac {\sqrt{\,2\,m}}{{\pi }} 
{\displaystyle \sum _{k=0}^{\infty }} \,{\displaystyle  {\frac{(-1
)^{k}}{\zeta^{\prime} ( - 2
\,j)}\,\Big [ {\displaystyle \frac {1}{j}}  + {\displaystyle \frac {2\,
\zeta^{\prime\,\prime} (- 2\,j)}{\zeta^{\prime} ( - 2\,j)}} \Big]}{}}  }} - {\displaystyle \frac {2\,
\sqrt{\,2\,m}}{\pi \,\zeta (2\,m)}} }
}
\label{Except_with_w=1}
\end{equation} 
where {\scriptsize ${\displaystyle { \sideset{}{\,^{'}}\sum_{n=1}}}$} means that the summation index $n$ excludes all instances in which $n=m(4h^2-1)$ and {\scriptsize ${\displaystyle \sum _{n=1}}^{\displaystyle \star}$} sums over all so-excluded values of $n$. In the final sum, the variable $j$ is defined by $j=m(2\,k+3)\,(2\,k+1)=m(4(k+1)^{2}-1)$ in terms of the summation index $k$.

\subsection{Case $a<0$}

In the case $a<0$, the bounding parabola of Figure \ref{fig:Fig1} opens to the left, and thereby $\mathfrak{S}$ always encloses all the poles of $1/\zeta(-2n)$. Therefore, the derivation conditions of (\ref{Master1_Zeta}) are not satisfied, and (\ref{Master1_Zeta}) must be modified by adding all the poles that, in this case reside in $\mathfrak{S}$, to the right-hand side (see Section 2). When this is done, \eqref{Master1_Zeta} becomes
\begin{equation}
\mapleinline{inert}{2d}{Int(w^(v^2)*w^(i*v)/Zeta(4*a*v^2+4*i*a*v)/cosh(Pi*v),v = -infinity ..
infinity) =
w^(1/4)/Zeta(a)+1/2*Pi*Sum(w^(-1/2*n/a)/Zeta(1,-2*n)/(-a^2-2*a*n)^(1/2
)/sinh(1/2*Pi*(-a^2-2*a*n)^(1/2)/a),n = ceil(-1/2*a) ..
infinity)-1/2*Pi*Sum(w^(-1/2*n/a)/Zeta(1,-2*n)/(a^2+2*a*n)^(1/2)/sin(1
/2*Pi*(a^2+2*a*n)^(1/2)/a),n = 1 .. floor(-1/2*a));}{%
\maplemultiline{
{\displaystyle \int _{ - \infty }^{\infty }} {\displaystyle 
\frac {w^{v^{2}+i\,v}}{\zeta (4\,a\,v\,(\,v + i))
\,\mathrm{cosh}(\pi \,v)}} \,dv=   
\mbox{}  {\displaystyle \frac {\pi}{2}} 
{\displaystyle \sum _{n=\mathrm{\lceil} - { {a
}/{2}} \rceil}^{\infty }} \,{\displaystyle \frac {w^{( - {n}/{2\,a
})}}{\zeta^{\prime} ( - 2\,n)\,\sqrt{ - a^{2} - 2\,a\,n}\,\mathrm{
sinh}({\displaystyle \frac {\pi \,\sqrt{ - a^{2} - 2\,a\,n}}{2\,a
}} )}}   \\
\mbox{} - {\displaystyle \frac {\pi}{2}} 
{\displaystyle \sum _{n=1}^{\mathrm{\lfloor} - { 
{a}/{2}}\rfloor}} \,{\displaystyle \frac {w^{( - {n}/{2\,a})
}}{\zeta^{\prime} (- 2\,n)\,\sqrt{a^{2} + 2\,a\,n}\,\mathrm{sin}(
{\displaystyle \frac {\pi \,\sqrt{a^{2} + 2\,a\,n}}{2\,a}} )}} 
{\,+\,\displaystyle \frac {w^{1/4}}{
\zeta (a)}}
  \hspace*{1.5cm} {a<0,\,a \neq -1,-2,\dots}\,\,,} 
\label{Srule2}
}
\end{equation}
where $\lfloor .. \rfloor$ and $\lceil.. \rceil$ signify the $floor$ and $ceiling$ operators respectively. In the case $a=-2\,m,m=1,2,\dots$,  the limiting variation becomes
\begin{equation}
\mapleinline{inert}{2d}{Int(w^(v^2)*w^(i*v)/Zeta(-8*m*v*(v+i))/cosh(Pi*v),v = -infinity ..
infinity) =
-1/8*w^(1/4)/m/Zeta(1,-2*m)*ln(w)+1/48*Pi^2*w^(1/4)/m/Zeta(1,-2*m)-1/2
*w^(1/4)*Zeta(2,-2*m)/Zeta(1,-2*m)^2-1/4*Pi/m^(1/2)*Sum(w^(1/4*n/m)*w^
(1/4)*w^(1/(4*m))/(n+1)^(1/2)/Zeta(1,-2*n-2*m-2)/sinh(1/2*Pi*(n+1)^(1/
2)/m^(1/2)),n = 0 ..
infinity)+1/4*Pi*Sum(w^(1/4*n/m)/Zeta(1,-2*n)/(m^2-m*n)^(1/2)/sin(1/2*
Pi*(m^2-m*n)^(1/2)/m),n = 1 .. m-1);}{%
\maplemultiline{
{\displaystyle \int _{ - \infty }^{\infty }} {\displaystyle 
\frac {w^{v^{2}+i\,v}}{\zeta ( - 8\,m\,v\,(v + i))\,
\mathrm{cosh}(\pi \,v)}} \,dv= - {\displaystyle \frac {1}{8}} \,
{\displaystyle \frac {w^{1/4}\,\mathrm{ln}(w)}{m\,\zeta^{\prime} ( - 2\,m)}}  + {\displaystyle \frac {1}{48}} \,{\displaystyle 
\frac {\pi ^{2}\,w^{1/4}}{m\,\zeta^{\prime} (- 2\,m)}}  - 
{\displaystyle \frac {1}{2}} \,{\displaystyle \frac {w^{1/4}\,
\zeta^{\prime\,\prime} (- 2\,m)}{\zeta^{\prime} (- 2\,m)^{2}}}  
\mbox{}\\ - {\displaystyle \frac {\pi}{4\,\sqrt{m}}} \,{\displaystyle  {
 {\displaystyle \sum _{n=0}^{\infty }} \,
{\displaystyle \frac {w^{\displaystyle ({n+1+m})/({4\,m})}}{\sqrt{n + 1}\,\zeta^{\prime} ( - 2\,n - 2\,m - 2)
\,\mathrm{sinh}({\displaystyle \frac {\pi \,\sqrt{1 + n }}{2\,
\sqrt{m}}} )}}  }{}}  
\mbox{}  \\+ {\displaystyle \frac {\pi}{4}}  
{\displaystyle \sum _{n=1}^{m - 1}} \,{\displaystyle \frac {w^{\displaystyle  {n}/{(4\,m)}}}{\zeta^{\prime} ( - 2\,n)\,\sqrt{m^{2} - m\,n}\,
\mathrm{sin}({\displaystyle  {\frac{\pi}{2} \,\sqrt{1 - n/m}}
} )}} \,.  }
}
\label{Srule2a}
\end{equation}
\newline 
If the integration contour in \eqref{Srule2a} is shifted upwards and parallel to the real $v$-axis, the left-opening parabola in the $z$-plane shifts to the right and opens rapidly, thereby enveloping any poles that are not already included. At $v \rightarrow v+i\infty$, contributions from the integral vanish, leaving a sum of contributions from the specified poles. After evaluating the corresponding residues, we are left with the sum rule
\begin{equation}
\mapleinline{inert}{2d}{Sum(1/2*exp(-1/8*tau*ln(w)/m)*2^(1/2)/Zeta(1,tau)/(2*m+tau)^(1/2)/sin
(1/4*Pi*2^(1/2)*(2*m+tau)^(1/2)/m^(1/2))+\{tau ->
conjugate(tau)\},tau)/m^(1/2) =
-1/24*w^(1/4)*Pi/m/Zeta(1,-2*m)+1/2*Sum(w^(1/4*(n+m+1)/m)/(n+1)^(1/2)/
m^(1/2)/Zeta(1,-2*n-2*m-2)/sinh(1/2*Pi*(n+1)^(1/2)/m^(1/2)),n = 0 ..
infinity)-1/2*Sum(w^(1/4*n/m)/Zeta(1,-2*n)/(m^2-m*n)^(1/2)/sin(1/2*Pi*
(m^2-m*n)^(1/2)/m),n = 1 ..
m-1)+1/4*1/Pi*w^(1/4)*ln(w)/m/Zeta(1,-2*m)+w^(1/4)*Zeta(2,-2*m)/Zeta(1
,-2*m)^2/Pi+4/Pi*Sum((-1)^k*w^(-1/4*(2*k+3)*(2*k+1))/Zeta(2*m*(2*k+3)*
(2*k+1)),k = 0 .. infinity);}{%
\maplemultiline{
{\displaystyle  {{\displaystyle \frac {1}{\sqrt{2\,m}}}{\displaystyle \sum _{\tau }}  
 \!  \,{\displaystyle \frac {\exp{( -\,{\displaystyle {\tau \,\mathrm{ln}(w)}/{8\,m})}}}{\zeta^{\prime} (\tau )\,\sqrt{2\,m + \tau }\,\mathrm{sin}({\displaystyle {\frac
{\pi}{2} \,\sqrt{1 + \tau/2m }}} )}}    }} = - 
{\displaystyle \frac {\pi}{24\,m}} \,{\displaystyle \frac {w^{1/4} }{\zeta^{\prime} (- 2\,m)}}  \\
\mbox{} + {\displaystyle \frac {1}{2\,\sqrt{m}}} 
{\displaystyle \sum _{n=0}^{\infty }} \,{\displaystyle \frac {w^{
 ({n + m + 1})/({4\,m})}}{\sqrt{n + 1}\,\zeta^{\prime} ( - 2\,n - 2\,m - 2)\,\mathrm{sinh}({\displaystyle \frac {\pi \,
\sqrt{n + 1}}{2\,\sqrt{m}}} )}}   \\
\mbox{} - {\displaystyle \frac {1}{2}}\,
{\displaystyle \sum _{n=1}^{m - 1}} \,{\displaystyle \frac {w^{ 
{n}/{4\,m}}}{\zeta^{\prime} (- 2\,n)\,\sqrt{m^{2} - m\,n}\,
\mathrm{sin}({\displaystyle {\frac{\pi}{2} \,\sqrt{1 - n/m}}{
}} )}}   + {\displaystyle \frac {1}{4}} \,
{\displaystyle \frac {w^{1/4}\,\mathrm{ln}(w)}{\pi \,m\,\zeta ^\prime(
- 2\,m)}}  
\mbox{} + {\displaystyle \frac {w^{1/4}\,\zeta^{\prime\prime} (- 2\,m)}{
\pi \, \zeta^{\prime} (- 2\,m)^{2}}} \\ + {\displaystyle \frac{4}{\pi}}\,{\displaystyle {
 {\displaystyle \sum _{k=0}^{\infty }} \,
{\displaystyle \frac {(-1)^{k}\,w^{\displaystyle - {(2\,k + 3)\,(2\,k
 + 1)}/{4}}}{\zeta (2\,m\,(2\,k + 3)\,(2\,k + 1))}} 
}{}} \,. }
}
\label{Srule_2a}
\end{equation}

\subsection{Case $a=i\beta, \beta>0$}

In this case, the parabola in the $z$-plane opens upwards and envelopes all zeros of $\zeta(z)$ at $z=\sigma+i\rho$ for all values of $\beta>1/(16\rho)$ (assuming $\sigma=1/2$).  However, all zeros that lead to poles belonging to $\Im(\tau)\leq 0$ are excluded. This leaves the following identity, valid for $w \leq 1$:
\begin{equation}
\mapleinline{inert}{2d}{Int(w^(v*(v+I))/Zeta(4*i*b*v*(v+I))/cosh(Pi*v),v = -infinity ..
infinity) =
w^(1/4)/Zeta(i*b)+1/2*Pi*Sum(1/(w^(-1/4*tau/i/b)*(-i*b+tau)^(1/2)*Zeta
(1,tau)*sinh(1/2*Pi*(-i*b*(i*b-tau))^(1/2)/i/b)),tau)/(i*b)^(1/2);}{%
\maplemultiline{
{\displaystyle \int _{ - \infty }^{\infty }} {\displaystyle 
\frac {w^{v\,(v + i)}}{\zeta (4\,i\,\beta\,v\,(v + i))\,\mathrm{
cosh}(\pi \,v)}} \,dv= 
{\displaystyle \frac {w^{1/4}}{\zeta (i\,\beta)}}  + 
{\displaystyle \frac {\pi}{2\,\sqrt{i\,\beta}}} \,{\displaystyle  {
 {\displaystyle \sum _{\substack {\Im(\tau)>\rho_{M} }}} \,{\displaystyle 
\frac {w^{ {\tau }/{4\,i\,\beta}}}{\sqrt{ \tau - i\,\beta }
\,\,\zeta^{\prime} (\tau )\,\mathrm{sinh}({\displaystyle  {\frac{\pi}{2} \,
\sqrt{-i\,\tau/\beta - 1 }}} )}} }{
}}  }
}
\label{ib_Int}
\end{equation}

In the case that $\beta<1/(16\rho_{M})$ (again, assuming $\sigma=1/2$) corresponding to the $M^{th}$ zero (in ordered  magnitude) of $\zeta(z)$ then all poles $\Im(\tau)<\Im{(\tau_{M})}$ must be  excluded from the sum. Numerically, for large values of $\rho_{M}$, this is almost impossible to detect because for correspondingly small values of $\beta$, the contribution of the sum over $\tau$ is negligible. Also, the utility of \eqref{ib_Int} is limited, since any attempt to translate the integration contour and transform the integration into a sum as done previously, fails because the series $S_{2}$ defined in \eqref{S2} does not converge unless $w>1$. Similar considerations apply to the case $\beta<0$.

\subsection{Case $a=\alpha+i\beta,\alpha<0,\beta>0$}

In the case of complex $a=\alpha+i\beta, \alpha<0, \beta>0$, the bounding parabola of Figure \ref{fig:Fig1}, which originally opened to the left for $\alpha<0, \beta=0$ (see above), rotates clockwise about the origin as $\beta$ increases, thereby truncating the infinite sum appearing in \eqref{Srule2}. At the same time, it quickly encloses poles belonging to non-trivial zeros of $\zeta(z)$, yielding the following result,
\begin{equation}
\mapleinline{inert}{2d}{Int(Pi/Zeta(4*a*v*(v+i))/cosh(Pi*v),v = -infinity .. infinity) =
Pi/Zeta(a)-1/2*Pi^2*Sum(1/(Zeta(1,-2*n)*(a^2+2*a*n)^(1/2)*sin(1/2*Pi*(
a^2+2*a*n)^(1/2)/a)),n = 1 ..
floor(-2*alpha*abs(a)^2/beta^2))-1/2*Pi^2*Sum(1/((-a+tau)^(1/2)*Zeta(1
,tau)*sinh(1/2*Pi*(-a*(a-tau))^(1/2)/a)),`Im(tau)>0`)/a^(1/2);}{%
\maplemultiline{
{\displaystyle \int _{ - \infty }^{\infty }} {\displaystyle 
\frac {1 }{\zeta (4\,a\,v\,(v + i))\,\mathrm{cosh}(\pi \,v)}} 
\,dv={\displaystyle \frac {1 }{\zeta (a)}}  
\mbox{} - {\displaystyle \frac {\pi }{2}}
{\displaystyle \sum _{n=1}^N} \,{\displaystyle \frac {1}{\zeta^{\prime} (- 2\,n)\,\sqrt{a^{
2} + 2\,a\,n}\,\mathrm{sin}\left ( {\displaystyle \frac {\pi \,\sqrt{a^{
2} + 2\,a\,n}}{2\,a}}\right ) }}   \\
\mbox{} - {\displaystyle \frac {\pi}{2\,\sqrt{a}}} \,{\displaystyle {
{\displaystyle \sum _{\Im(\tau_{M_{-}}) \leq \Im(\tau)\leq \Im(\tau_{M_{+}}) }} 
\,{\displaystyle \frac {1}{\sqrt{ - a + \tau }\,\zeta^{\prime} (\tau 
)\,\mathrm{sinh}\left({\displaystyle \frac {\pi \,\sqrt{ - a\,(a - 
\tau )}}{2\,a}} \right)}} }{}}  }
}
\label{nega_cplx}
\end{equation}
where 
\begin{equation}
\mapleinline{inert}{2d}{beta = -(-(4*alpha+2*n)*alpha)^(1/2)*alpha/(2*alpha+n);}{%
\[
\beta \neq - \alpha\,{\displaystyle \frac {\sqrt{ - (4\,\alpha  + 2\,L)\,
\alpha }\, }{2\,\alpha  + L}}, 
\]
}
\label{BetaNeq}
\end{equation}
 $L$ being a positive integer,
\begin{equation}
N={{\lfloor} - {\displaystyle 
\frac {2\,\alpha \, \left|  \! \,a\, \!  \right| ^{2}}{\beta ^{2}
}}\rfloor}\,,
\end{equation}

and the sum extends over values of $\tau$ labelling $\zeta(\tau)=0$ that are interior to the parabola for particular values of $\alpha$ and $\beta$ bounded by $\Im(\tau_{M_{\pm}})$ - see below. In the case that \eqref{BetaNeq} is false, the integral diverges and must be interpreted in the sense of a principal value. In that case, \eqref{nega_cplx} reduces to an identity between the residue of the integral on the left, and a corresponding term on the right belonging to the residue of the pole through which the contour of integration passes, related via the functional equation. For example, if $\alpha=-1/2,\, \beta=-1/2, L=2$, the exceptional case reduces to
\begin{equation}
\mapleinline{inert}{2d}{Zeta(5) = 4/3*Pi^4*Zeta(1,-4);}{%
\[
\zeta (5)={\displaystyle \frac {4}{3}} \,\pi ^{4}\,\zeta^{\prime}(-4
)
\]
}
\label{zeta5}
\end{equation}
which can be independently obtained by once-differentiating the functional equation for $\zeta(z)$ and evaluating the limit $z\rightarrow5$.
\newline

 From Section 2, and assuming that all zeros lie on the critical line \cite{Odlyzko}, these limiting values of the $\tau$ sum in \eqref{nega_cplx} are bounded by
\begin{equation}
\mapleinline{inert}{2d}{Lim :=
1/2*beta*(4*abs(a)^2+alpha)/alpha^2+abs(a)^2*(4*beta^2+2*alpha)^(1/2)/
alpha^2;}{%
\[
\tau_{M_{\pm}}= {\displaystyle \frac {1}{2}} \,{\displaystyle 
\frac {\beta \,(4\, \left|  \! \,a\, \!  \right| ^{2} + \alpha )
}{\alpha ^{2}}}  \pm {\displaystyle \frac { \left|  \! \,a\, \! 
 \right| ^{2}\,\sqrt{4\,\beta ^{2} + 2\,\alpha }}{\alpha ^{2}}} 
\]
}
\label{taupm}
\end{equation}

As an example, for particular, carefully chosen values of $\alpha$ and $\beta$, both sums in \eqref{nega_cplx} will contain only one term, leading to interesting results such as, for the case $\alpha=-1/2, \beta=1$:
\begin{equation}
\mapleinline{inert}{2d}{int(1/(Zeta((-2+4*i)*v*(v+i))*cosh(Pi*v)),v = -infinity .. infinity)
=
1/Zeta(-1/2+i)-4*Pi^3/Zeta(3)/(-7+4*i)^(1/2)/sin((1/10+1/5*i)*Pi*(-7+4
*i)^(1/2))+2*Pi/(2-4*i+4*tau[1])^(1/2)/Zeta(1,tau[1])/sinh((1/5+2/5*i)
*Pi*((1/2-i)*(-1/2+i-tau[1]))^(1/2))/(-2+4*i)^(1/2);}{%
\maplemultiline{
{\displaystyle \int _{ - \infty }^{\infty }} {\displaystyle 
\frac {1}{\zeta (( - 2 + 4\,i)\,v\,(v + i))\,\mathrm{cosh}(\pi \,
v)}} \,dv={\displaystyle \frac {1}{\zeta ( - {\displaystyle 
 {1}/{2}}  + i)}}  
\mbox{} - {\displaystyle \frac {4\,\pi ^{3}}{\zeta (3)\,\sqrt{ - 
7 + 4\,i}\,\mathrm{sin} \left( ({\displaystyle {1}/{10}}  + 
{\displaystyle  {i}/{5}} )\,\pi \,\sqrt{ - 7 + 4\,i} \right ) }}  \\
\mbox{} + {\displaystyle \frac {2\,\pi }{\sqrt{2 - 4\,i + 4\,{
\tau _{1}}}\,\zeta^{\prime} ({\tau _{1}})\,\mathrm{sinh}\left ( (
{\displaystyle {1}/{5}}  + {\displaystyle {2\,i}/{5}} )
\,\pi \,\sqrt{({\displaystyle {1}/{2}}  - i)\,( - 
{\displaystyle {1}/{2}}  + i - {\tau _{1}})}\right )\,\sqrt{ - 2 + 
4\,i}}}  }
}
\label{amhalf}
\end{equation} 
where we have used the identity \eqref{ZetaprimeNegN}. We note that most of the numerical value of the integrand is contained in the range $0<v\leq1$ for this combination of $\alpha$ and $\beta$ which corresponds to the second quadrant of the $z-$plane bounded by $-7<\Re(z)<0$ and $0<\Im(z)<2$. Thus the general result \eqref{nega_cplx} defines a numerical relationship between numerically dominant values of $\zeta(z)$ in one part of the complex $z$-plane and zeros that appear elsewhere, without recourse to the functional equation. As previously noted, it is possible to shift the integration contour parallel to the real $v$ axis upwards to infinity where the integral vanishes, adding terms corresponding to the residues of poles encountered during transit. The result is the same as \eqref{Srule1}, except that the $\tau$ sum extends over all zeros of $\zeta(\tau)$ including $\Im(\tau)<0$.

\subsection{General case $a=\alpha+i\beta$}

The above considerations lead to interest in the general case $a= \alpha+i\beta$, where, for the remainder of this section, unless stated otherwise, we take $\alpha,\beta>0$. In this case, the illustrative parabola shown in Fig. \ref{fig:Fig1} will be skewed to the right or left, but will open upward, and, depending on the relative magnitude of $\alpha$ and $\beta$, it may or may not enclose poles of $1/\zeta(z)$. In general, the result \eqref{ib_Int} will apply, with the replacement $i\beta\rightarrow \alpha+i\beta$, and the value of the parameter $M$ which labels those poles that are enclosed by the parabola must be determined by reference to Section 2. In the following sections several examples will be given. Since the parameter $\alpha$ does not vanish, when the integration contour of \eqref{ib_Int} is translated upwards to infinity in the $v$-plane, the resulting series \eqref{S2} now converges when $w=1$. In addition, this translation will envelop upper half-plane poles that were previously not enclosed, so that when the residues are calculated, the resulting expressions will contain a sum over all poles belonging to zeros of $\zeta(z)$ with the exception of those that were previously omitted. Effectively this results in a term that cancels the extra (limited) sum that appears in \eqref{ib_Int}, and the sum rule is the same as \eqref{Srule1}, with the proviso that $a$ is replaced by $\alpha+i\beta$. In the case $\alpha>0, \beta<0$ similar comments apply, except that the parabola opens downward and to the right. In the end, the general result \eqref{Srule1} applies; \eqref{ib_Int} also follows except that the summation constraint $\Im(\tau)>0$ becomes $\Im(\tau)<\rho_{M}$. 

\section{Examples}

In the following, we shall apply the principles enunciated above to a number of different choices of function $F(z)$. To maintain clarity, the left-hand side of all results represents contributions from singularities encountered when the integration contour is translated to $+i\,\infty$ parallel to the real $v-$axis, except where they may add to or cancel contributions from poles already included in the right-hand side. Such occurrences   will be discussed in the text as they may appear. The right-hand side represents singularities contributed by poles that originally resided in $\mathfrak{S}$.
\subsection{Case $a=1/2+i\rho_{m}$}

With the knowledge that many points corresponding to $\zeta(\tau_{m})=0$ lie on the critical line \cite{Odlyzko}, consider the case that the parameter $a\approx \tau_m$ where $ \tau_m=1/2+i\rho_{m}$, the $m^{th}$ non-trivial zero of $\zeta(z)$. First of all, we find that, for this case, a large number of  poles $1/\zeta(\tau)$ referenced in \eqref{Master1_Zeta} are enclosed by the parabolic region in the $z$ plane (see Section 2), so an additional sum corresponding to the residues of such terms (residing in $\mathfrak{S}$), must be added to \eqref{Master1_Zeta} (similar to \eqref{Srule2}). Specifically, and to maintain simplicity by letting $w=1$, we have, with $a\approx \tau_m$
\begin{equation}
\mapleinline{inert}{2d}{Int(1/(Zeta(4*a*v*(v+i))*cosh(Pi*v)),v = -infinity .. infinity) =
-1/24*(Pi^2*Zeta(1,a)+12*Zeta(2,a)*a)/a/Zeta(1,a)^2+1/2*Pi*Sum(1/((-a*
(a-tau))^(1/2)*Zeta(1,tau)*sinh(1/2*Pi*(-a*(a-tau))^(1/2)/a)),`Im(tau)
>0`);}{%
\maplemultiline{
{\displaystyle \int _{ - \infty }^{\infty }} {\displaystyle 
\frac {1}{\zeta (4\,a\,v\,(v + i))\,\mathrm{cosh}(\pi \,v)}} \,dv
= 1/\zeta(a) 
\mbox{} + {\displaystyle \frac {\pi}{2}} 
{\displaystyle \sum^{\tau_{_{M}}} _{\substack { {\mathit{\Im(\tau)>0}}}}} \,{\displaystyle 
\frac {1}{\sqrt{ - a\,(a - \tau )}\,\zeta^{\prime} (\tau )\,\mathrm{
sinh}({\displaystyle \frac {\pi \,\sqrt{ - a\,(a - \tau )}}{2\,a}
} )}}   }
}  
\label{a_approx_tau}
\end{equation}   
where the summation index $\tau=\sigma+i\rho$ extends over all zeros of $\zeta(\tau)$ such that $\rho>0$ and $\Im(\tau_{M})<\Im(z)$. In the case of equality (viz. $a=1/2+i\rho_{1}$), with $\tau_{M}=1/2+i\rho_{M}$ it turns out that $\rho_{M}< 45,268$, defined by the intersection of the integration contour in the $z$-plane and the line $\Re(z)=\frac{1}{2}$, corresponding to $M=56,791$. Notice that the bounding curves of the parabolic region in the $z$-plane enclose the poles in the opposite direction to that in which poles are enclosed in $\mathfrak{S}$, hence the sign of the summation term in \eqref{a_approx_tau}. In the general case of equality, that is $a_{m}\equiv \tau_{m} =\sigma_{m}+i\rho_{m}$ a limiting cancellation exists between the first term on the right of   \eqref{a_approx_tau} and the corresponding term of the sum belonging to $\tau=\tau_{m}$. The result, again setting $w=1$ for simplicity is

\begin{equation}
\mapleinline{inert}{2d}{Int(1/(Zeta(4*a*v*(v+i))*cosh(Pi*v)),v = -infinity .. infinity) =
-1/24*(Pi^2*Zeta(1,a)+12*Zeta(2,a)*a)/a/Zeta(1,a)^2+1/2*Pi*Sum(1/((-a*
(a-tau))^(1/2)*Zeta(1,tau)*sinh(1/2*Pi*(-a*(a-tau))^(1/2)/a)),`Im(tau)
>0`);}{%
\maplemultiline{
{\displaystyle \int _{ - \infty }^{\infty }} {\displaystyle 
\frac {1}{\zeta (4\,\tau_{m}\,v\,(v + i))\,\mathrm{cosh}(\pi \,v)}} \,dv
= - {\displaystyle \frac {1}{24}} \,{\displaystyle \frac {\pi ^{2
}\,\zeta^{\prime} (\tau_{m}) + 12\,\zeta^{\prime\prime} (\tau_{m})\,\tau_{m}}{\tau_{m}\,\zeta^{\prime} (\tau_{m})^{2}}
}  \\
\mbox{} + {\displaystyle \frac {\pi}{2}}   
{\displaystyle \sum^{\tau_{_{M}}} _{\substack {{\mathit{\Im(\tau)>0}}}}}^{_{\displaystyle _\prime}} \,{\displaystyle 
\frac {1}{\sqrt{ - \tau_{m}\,(\tau_{m} - \tau )}\,\zeta^{\prime} (\tau )\,\mathrm{
sinh}({\displaystyle \frac {\pi \,\sqrt{ - \tau_{m}\,(\tau_{m} - \tau )}}{2\,\tau_{m}}
} )}}  }
}
\label{a=tau}
\end{equation}

where the summation omits the index $\tau=\tau_{m}$, indicated by the $^{\prime}$ symbol. Seen another way, the limiting cancellation of divergent terms to yield \eqref{a=tau} corresponds to the coalescence of two poles, one belonging to $1/\cosh(\pi\,v)$, the second to $1/\zeta(z)$,  generating a dipole at $v=-i/2$ inside $\mathfrak{S}$, when $a=\tau_{m}$.
\newline

As in the previous sections, the integral on the left-hand side of \eqref{a=tau} can be moved upwards in the complex $v$-plane, and the residues of the poles encountered as it does so must be added to the left-hand side. It turns out that the corresponding contour in the $z$-plane now encloses poles belonging to complex conjugate zeros $\zeta{(\overline{\tau})}$ as well as the corresponding sums that have been seen previously. Hence a new sum including such terms must be added to \eqref{a=tau}. Of interest, a lack of symmetry means that the limiting index of the sums belonging to $\zeta({\tau})$ and $\zeta(\bar{\tau})$ will differ, as the contour moves up the imaginary v-axis. For example, for the case $a=\tau_{1}$, when the contour reaches $v=15i$ the upper limit M in \eqref{a=tau} is approximately $7\times 10^9$ whereas the upper limit in the corresponding sum containing residues belonging to $\zeta(\bar{\tau})$ is 14,438. Of course, when the countour has been translated to infinity, the integral vanishes, and the corresponding sums develop infinite limits. The final result is
\begin{equation}
\mapleinline{inert}{2d}{Sum(1/((-tau[m]+tau)^(1/2)*Zeta(1,tau)*sinh(1/2*Pi*(-tau[m]*(tau[m]-t
au))^(1/2)/tau[m]))+1/((-tau[m]+conjugate(tau))^(1/2)*Zeta(1,conjugate
(tau))*sinh(1/2*Pi*(-tau[m]*(tau[m]-conjugate(tau)))^(1/2)/tau[m])),`I
m(tau[M])>0`) = 4/Pi*Sum((-1)^k/Zeta(-tau[m]*(2*k+3)*(2*k+1)),k = 0 ..
infinity)*tau[m]^(1/2)+Sum(1/((tau[m]+2*n)^(1/2)*Zeta(1,-2*n)*sin(1/2*
Pi*(tau[m]+2*n)^(1/2)/tau[m]^(1/2))),n = 1 ..
infinity)-1/((-tau[m]+conjugate(tau[m]))^(1/2)*Zeta(1,conjugate(tau[m]
))*sinh(1/2*Pi*(-tau[m]*(tau[m]-conjugate(tau[m])))^(1/2)/tau[m]))+1/1
2*1/Zeta(1,tau[m])*Pi/tau[m]^(1/2)+tau[m]^(1/2)*Zeta(2,tau[m])/Pi/Zeta
(1,tau[m])^2;}{%
\maplemultiline{
{\displaystyle \sum ^{\infty} _{\mathit{\Im(\tau)>0}}}^{_{\displaystyle _\prime}}  \left( {\vrule 
height1.50em width0em depth1.50em} {\displaystyle 
\frac {1}{\sqrt{ - {\tau _{m}} + \tau }\,\zeta^{\prime}(\tau )\,
\mathrm{sinh}({\displaystyle \frac {1}{2}} \,{\displaystyle 
\frac {\pi \,\sqrt{ - {\tau _{m}}\,({\tau _{m}} - \tau )}}{{\tau 
_{m}}}} )}}  \mbox{} +\{\tau\rightarrow \overline{\tau}\} \right ) \\ = 
{\displaystyle  {\frac{4\,\sqrt{{\tau _{m}}}} {\pi} {\displaystyle \sum _{k=0}^{
\infty }} \,{\displaystyle \frac {(-1)^{k}}{\zeta ( - {\tau _{m}}
\,(2\,k + 3)\,(2\,k + 1))}}  }{
}}  
\mbox{} + {\displaystyle \sum _{n=1}^{\infty }} \,
{\displaystyle \frac {1}{\sqrt{{\tau _{m}} + 2\,n}\,\zeta^{\prime}(
 - 2\,n)\,\mathrm{sin}({\displaystyle \frac {1}{2}} \,
{\displaystyle \frac {\pi \,\sqrt{{\tau _{m}} + 2\,n}}{\sqrt{{
\tau _{m}}}}} )}}  \\
\mbox{} - {\displaystyle \frac {1}{\sqrt{ - {\tau _{m}} + 
\overline{({\tau _{m}})}}\,\zeta^{\prime}(\overline{({\tau _{m}})})
\,\mathrm{sinh}({\displaystyle \frac {1}{2}} \,{\displaystyle 
\frac {\pi \,\sqrt{ - {\tau _{m}}\,({\tau _{m}} - \overline{{
\tau _{m}})}}}{{\tau _{m}}}} )}}  + {\displaystyle \frac {1}{12}
} \,{\displaystyle \frac {\pi }{\zeta^{\prime}({\tau _{m}})\,\sqrt{{
\tau _{m}}}}} 
\mbox{} + {\displaystyle \frac {\sqrt{{\tau _{m}}}\,\zeta (2, \,{
\tau _{m}})}{\pi \,\zeta^{\prime}({\tau _{m}})^{2}}}  }
}
\label{a=taum}
\end{equation} 
A few comments are relevant to \eqref{a=taum}. The sum on the left-hand side extends over all indices $\tau$ labelling $\zeta(\tau)=0$ except for the element $\tau=\tau_{m}$. In that case, the singularity belonging to $\Im(\tau_{m})>0$ is cancelled as in \eqref{a=tau} and the cancellation terms have been written explicitly on the right. Since there is no singularity belonging to $\Im({\overline{\tau_{m}})} \equiv \Im(\tau_{m})<0$, that term has been removed from the sum and also explicitly included on the right-hand side of \eqref{a=taum}. In that manner the sum on  the left-hand side can be written symmetrically, omitting only the index $\tau=\tau_{m}$, indicated by the $^{\prime}$ symbol. 
\newline

A numerical study of this sum is interesting. For large values of $\rho$, it can be shown that 
\begin{equation}
\mapleinline{inert}{2d}{T(tau[m],tau) :=
1/((-tau[m]+tau)^(1/2)*sinh(1/2*Pi*(-tau[m]*(tau[m]-tau))^(1/2)/tau[m]
)) =
(1-i)*2^(1/2)*exp(-c[1]*Pi*rho^(1/2))*exp(i*c[2]*Pi*rho^(1/2))/rho^(1/
2);}{%
\[
\mathrm{\widehat{T}}({\tau _{m}}, \,\tau ) \equiv {\displaystyle \frac {1}{
\sqrt{ - {\tau _{m}} + \tau }\,\mathrm{sinh}({\displaystyle 
\frac {1}{2}} \,{\displaystyle \frac {\pi \,\sqrt{ - {\tau _{m}}
\,({\tau _{m}} - \tau )}}{{\tau _{m}}}} )}} \approx{\displaystyle 
\frac {\sqrt{2}\,(1 - i)\,\exp{( - ({c_{1}-i\,c_{2})}\,\pi \,\sqrt{\rho })}}{\sqrt{\rho }}} 
\]
}
\label{tp}
\end{equation}
and
\begin{equation}
\mapleinline{inert}{2d}{T(tau[m],tau) :=
1/((-tau[m]+tau)^(1/2)*sinh(1/2*Pi*(-tau[m]*(tau[m]-tau))^(1/2)/tau[m]
)) =
(1-i)*2^(1/2)*exp(-c[1]*Pi*rho^(1/2))*exp(i*c[2]*Pi*rho^(1/2))/rho^(1/
2);}{%
\[
\mathrm{\widehat{T}}({\tau _{m}}, \,\overline{\tau} ) \approx{\displaystyle 
\frac {\sqrt{2}\,(1 + i)\,\exp{( ({c_{2}+i\,c_{1})}\,\pi \,\sqrt{\rho })}}{\sqrt{\rho }}} 
\]
}
\label{tm}
\end{equation}
where
\begin{eqnarray}
{c_{1}}={\displaystyle \frac {1}{2}} \,{\displaystyle \frac {1 + 
2\,{z_{2}}^{2}\,{\rho_{m}}}{{z_{2}}\,( - {z_{2}}^{2} + 2\,{\rho_{m}})^{
2}}} \nonumber
\\
{c_{2}}={\displaystyle \frac {1}{2}} \,{\displaystyle \frac {1}{(
 - {z_{2}}^{2} + 2\,{\rho_{m}})\,{z_{2}}}}
 \\
 {z_{2}}=\sqrt{\sqrt{4\,{\rho_{m}}^{2} + 1} + 2\,{\rho_{m}}} \nonumber
 \label{c1c2z2}
\end{eqnarray}
 
Using  $\rho_{m}=\rho_{1}=14.134725142$, we find \ignore{$c_{1}=0.4176115889$} $c_{1}=0.418$ and \ignore{$c_{2}=-0.007383960883$} $c_{2}=-0.00738$, which leads to reasonably quick convergence of the first terms inside the summation on the left-hand side of \eqref{a=taum}. However, the conjugate term in this same sum is extraordinarily slow to converge. \ignore{Figure 2 illustrates the magnitude of the real and imaginary parts of $\mathrm{\widehat{T}}({\tau _{m}}, \,\overline{\tau} )$ over a reasonable(?!) range of $\rho$ using $\tau_{m}=\tau_{1}$. For larger choices of $\tau_{m}$, convergence is more lethargic.} Note that these observations do not include the effect of the denominator function $\zeta^{\prime}(\tau)$ in  \eqref{a=taum} whose absolute value more-or-less increases slowly with $\rho$. Under the assumption that the zeros labelled by $\tau$ are to be found only on the critical line \cite{Odlyzko}, we have investigated the veracity of \eqref{a=taum} for the case $m=1$ courtesy of the Mathematica \cite{Math} supplied database $\textbf{ZetaZero[...]}$, finding, with default arithmetic precision, that the series on the left shows signs of convergence when 400,000 terms are included in the sum (see Figure \ref{fig:Fig3}). Similarly, if the integral on the left-hand side of \eqref{a_approx_tau} is evaluated numerically for small values of upwards translation in the complex $v$-plane, the sums labelled by $k$ and $n$ indices in \eqref{a=taum} are also truncated because the bounding contour in the $z$ plane cuts the $\Im(v)=0$ axis. This has been verified numerically for a small number of upward translations with $0\leq \Im(v)\leq 15i$. 
\begin{figure}
     \caption{Fractional error in the Real (dotted) and Imaginary (solid) parts of the first 400,000 terms of the partial sum defined by the left-hand side of \eqref{a=taum} for the case $m=1$.}
      \label{fig:Fig3}
       \includegraphics[keepaspectratio,totalheight=4.5in]{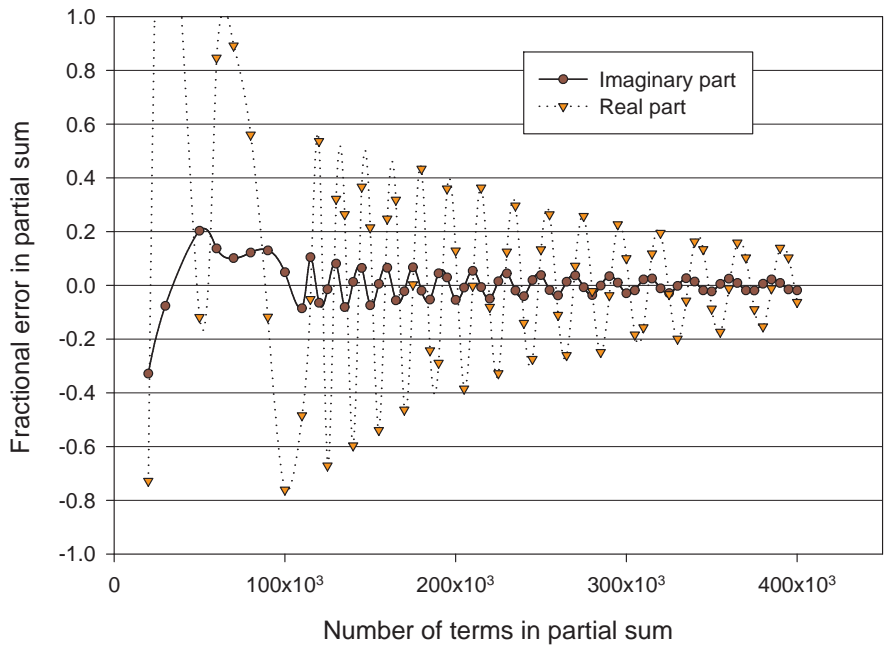}
     \end{figure}

\subsection{Interesting Variations}

A number of interesting variations immediately arise from the foregoing. Consider the case $a=1/2$, and re-evaluate the left hand side of \eqref{Srule1} by breaking it into its real and imaginary parts. We shall use the symbol $``\doteq "$ in the following to denote those equalities that are only valid at a zero on the critical line.  Furthermore, for the remainder of this subsection, we shall use $\tau=1/2+i\rho$ and refer to notations \cite{ReynaDeLune} 
\begin{equation}
\zeta (\tau)= e^{-i\phi(\rho)}Z(\rho) 
\label{Zeta_phase}
\end{equation}
and \cite{Milgram} 
\begin{equation}
\zeta (\tau)= \zeta_{R}(\tau) + i \zeta_{I}(\tau)
\end{equation}
giving
\begin{equation}
\phi(\rho)=-\arctan(\zeta_{I}(\tau)/\zeta_{R}(\tau)) \pm k\pi.
\end{equation}
where $k=0,1,...$ and the equivalence between the two representations can be verified by straightforward, but tedious calculation.
It has been shown \cite[Eq.(147)]{Milgram} that, if $\zeta(\tau)=0$, then
\begin{equation}
\frac{\Re(\zeta^{\prime}(\tau))}{\Im(\zeta^{\prime}(\tau))}\doteq\frac{N(\rho)}{D_{R}(\rho)}\,.
\end{equation}
where $N(\rho)$, $D_{R}(\rho)$ and various related symbols $D_{+}$ and $D_{-}$ - see below- are defined in the Appendix. Define $T(a,\tau)$, a representative single term in the sum over zeros appearing in \eqref{Srule1}. 
\begin{equation}
{T(a,\tau)\equiv \displaystyle } \,{\displaystyle \frac {1}{\sqrt{ - 
a + \tau }\,\zeta^{\prime} (\tau )\,\mathrm{sinh}({\displaystyle 
\frac {\pi \,\sqrt{ - a\,(a - \tau )}}{2\,a}} )}} +\{\tau \rightarrow \overline{(\tau)}\}
\label{Tdef}
\end{equation}

After arduous calculation, using Maple \cite{Maple}, we arrive at

\begin{equation}
\mapleinline{inert}{2d}{1/2*1/((-a+tau)^(1/2)*Zeta(1,tau)*sinh(1/2*Pi*(-a*(a-tau))^(1/2)/a))
=
2^(1/2)*Re(Zeta(1,tau))*((U(rho)-1)*cos(1/2*Pi*rho^(1/2))*sinh(1/2*Pi*
rho^(1/2))+(U(rho)+1)*sin(1/2*Pi*rho^(1/2))*cosh(1/2*Pi*rho^(1/2)))/rh
o^(1/2)/(cosh(Pi*rho^(1/2))-cos(Pi*rho^(1/2)))/abs(Zeta(1,tau))^2;}
{%
\maplemultiline{
T(1/2,\tau)     \doteq 
{\displaystyle \frac {2\,\sqrt{2}\,\Re (\zeta^{\prime} (\tau ))\,((
\mathrm{U}(\rho ) - 1)\,\mathrm{cos}({\displaystyle \frac {\pi \,
\sqrt{\rho }}{2}} )\,\mathrm{sinh}({\displaystyle \frac {\pi \,
\sqrt{\rho }}{2}} ) + (\mathrm{U}(\rho ) + 1)\,\mathrm{sin}(
{\displaystyle \frac {\pi \,\sqrt{\rho }}{2}} )\,\mathrm{cosh}(
{\displaystyle \frac {\pi \,\sqrt{\rho }}{2}} ))}{\sqrt{\rho }\,(
\mathrm{cosh}(\pi \,\sqrt{\rho }) - \mathrm{cos}(\pi \,\sqrt{\rho
 }))\, \left|  \! \,\zeta^{\prime} (\tau )\, \!  \right| ^{2}}} }
}
\label{Q9}
\end{equation}
where
\begin{equation}
\mapleinline{inert}{2d}{U(rho) =
cot(-1/2*argument(GAMMA(1/2+i*rho))+1/2*rho*ln(2*Pi)+1/2*arctan(tanh(1
/2*Pi*rho)));}{%
\[
\mathrm{U}(\rho )=\mathrm{cot}( - {\displaystyle \frac {1}{2}} \,
\mathrm{argument}(\Gamma ({\displaystyle \frac {1}{2}}  + i\,\rho
 )) + {\displaystyle \frac {1}{2}} \,\rho \,\mathrm{ln}(2\,\pi )
 + {\displaystyle \frac {1}{2}} \,\mathrm{arctan}(\mathrm{tanh}(
{\displaystyle \frac {\pi \,\rho }{2}} )))\,\,.
\]
}
\label{U(rho)}
\end{equation}
Consistently, reiterating that we are using the common assumption that $\zeta^{\prime}(z)$ does not vanish anywhere on the critical line, further simplification is available by noting that \cite{Milgram}
\begin{equation}
\mapleinline{inert}{2d}{Re(Zeta(1,1/2+rho*I))/abs(Zeta(1,1/2+rho*I))^2 =
D[r]/Re(Zeta(1,1/2+rho*I));}{%
\[
{\displaystyle \frac {\Re (\zeta^{\prime} (\tau))}{ \left|  \! \,\zeta^{\prime} (\tau)\, \!  \right| ^{2}}} \doteq{\displaystyle 
\frac {{\mathrm{D}_{R}(\rho)}}{\Re (\zeta^{\prime} (\tau))}} 
\]
}
\label{Re/abs(Zeta)}
\end{equation}
which reduces the dependence of the sum rule from the properties of $\zeta^{\prime}(\tau)$ to that of its real part only, giving

\begin{equation}
\mapleinline{inert}{2d}{Q5Pr :=
1/4*2^(1/2)*(`D[+]`*cos(1/2*Pi*rho^(1/2))*sinh(1/2*Pi*rho^(1/2))+`D[-]
`*sin(1/2*Pi*rho^(1/2))*cosh(1/2*Pi*rho^(1/2)))/(-cosh(Pi*rho^(1/2))+c
os(Pi*rho^(1/2)))/rho^(1/2)/Pi^(1/2)/Re(Zeta(1,1/2+i*rho));}{%
\[
{T(1/2,\tau) \doteq \displaystyle \frac {\sqrt{2}}{2\,\sqrt{\pi\,\rho }}} \,{\displaystyle 
\frac {(\mathit{D_{+}}\,\mathrm{cos}({\displaystyle 
\frac {\pi \,\sqrt{\rho }}{2}} )\,\mathrm{sinh}({\displaystyle 
\frac {\pi \,\sqrt{\rho }}{2}} ) + \mathit{D_{-}}\,\mathrm{sin}(
{\displaystyle \frac {\pi \,\sqrt{\rho }}{2}} )\,\mathrm{cosh}(
{\displaystyle \frac {\pi \,\sqrt{\rho }}{2}} ))}{( - \mathrm{
cosh}(\pi \,\sqrt{\rho }\,) + \mathrm{cos}(\pi \,\sqrt{\rho }\,))\;
\Re\, (\zeta^{\prime} ({\displaystyle 
 {1}/{2}}  + i\,\rho ))}} 
\]
}
\label{Q5}
\end{equation}

Although \eqref{Q9} is valid only at $\zeta(\tau)=0$, the right-hand side is a valid function of $\rho$ in it's own right, leading to the possibility that $\rho$ satisfies
\begin{equation}
 - {\displaystyle \frac {1}{2}} \,
\mathit{arg}(\Gamma ({\displaystyle \frac {1}{2}}  + i\,\rho
 )) + {\displaystyle \frac {1}{2}} \,\rho \,\mathrm{ln}(2\,\pi )
 + {\displaystyle \frac {1}{2}} \,\mathrm{arctan}(\mathrm{tanh}(
{\displaystyle \frac {\pi \,\rho }{2}}))=(n+1/2)\pi\,\,,
\label{U(rho)=inf}
\end{equation}
corresponding to the case where $U(\rho)$ in \eqref{U(rho)} diverges.
\newline

Any solution $\rho_{s}$ of \eqref{U(rho)=inf} that also satisfies $\zeta(1/2+i\rho_{s})=0$ , would result in an inconsistency in \eqref{Srule1}, since the right-hand side of that relation is finite for $a=1/2$ but (one of the terms on) the left-hand side would diverge. In fact, the equivalence of \eqref{Q5} and \eqref{Q9} demonstrates that the factor $\Re (\zeta^{\prime}(\tau))$ in \eqref{Q9} (or $D_{R}(\rho)$ in \eqref{Re/abs(Zeta)}) must vanish when \eqref{U(rho)=inf} is satisfied, thereby cancelling a potential divergence. Equivalently, since $D_{R}$ carries the $half-zero$ (defined by $\Re(\zeta^{\prime}(1/2+i\rho))=0 \; ;\; \Im(\zeta^{\prime}(1/2+i\rho))\neq 0)$ - see \eqref{Re/abs(Zeta)}), it can be shown that the only simultaneous solution when both $D_{R}=0$ and the numerator terms of \eqref{Q5} vanish, requires that $\Re(\Gamma(1/2+i\rho))=1$, whose only solution in turn occurs when $\rho=0.4102$, which does not correspond to a zero of $\zeta(\tau)$. 
\newline

The significance of this observation is that any point  $\rho_{0}$ satisfying the half-zero condition $\Re(\zeta^{\prime}(1/2+i\rho_{0}))=0$ will never be found to also satisfy $\zeta(1/2+i\rho_{0})=0$, since such a coincidence - an uncancelled divergence of \eqref{Q5} - would lead to a contradiction of \eqref{Q9}. Expressed another way, as a function of $\rho$, since

\begin{equation}
\zeta^{\prime}(1/2+i\rho)\doteq e^{-i\phi(\rho)}Z(\rho) = \zeta_{R}^{\prime}(\tau) + i \zeta_{I}^{\prime}(\tau)
\label{Zeta_phase_prime}
\end{equation} 

the phase $\phi(\rho)$ in \eqref{Zeta_phase} or \eqref{Zeta_phase_prime} will never assume the value(s) $\pm\pi/2$ at a point $\rho$ satisfying $\zeta(1/2+i\rho)=0$. This is consistent with our use of the hypothesis that the zeros of $\zeta(1/2+i\rho)$ are simple and, self-consistently, suggests the negative converse of Titchmarsh's proposition quoted in \cite[Proposition 10]{ReynaDeLune}: "If $\zeta^{\prime}(1/2+i\rho) = 0$, then $\zeta(1/2+i\rho) = 0$."  That is: {\it{If $\zeta(1/2+i\rho) = 0$, then $\zeta^{\prime}(1/2+i\rho) \neq 0$}}. 
\newline

Finally, considering the asymptotic properties of \eqref{Q9}, it can be shown that for $\rho \rightarrow \infty$, 
\begin{equation}
\mapleinline{inert}{2d}{T(1/2,1/2+i*rho) =
1/4*(-sin(-rho*ln(4*Pi/(1+4*rho^2)^(1/2))-1/2*Pi*rho^(1/2)-rho+1/2*arc
tan(2*rho))+cos(-rho*ln(4*Pi/(1+4*rho^2)^(1/2))-1/2*Pi*rho^(1/2)-rho+1
/2*arctan(2*rho))-sin(1/2*Pi*rho^(1/2))+cos(1/2*Pi*rho^(1/2)))*exp(-1/
2*Pi*rho^(1/2))*2^(1/2)/rho^(1/2)/Re(Zeta(1,1/2+i*rho));}{%
\maplemultiline{ \lim\limits_{\rho \to \infty}
\mathrm{T}({ \frac{1}{2}} , \,{ 
\frac{1}{2}}  + i\,\rho ) \doteq {\displaystyle \frac {\sqrt{2}\exp{( - {\displaystyle  {\pi \,\sqrt{\rho }}/{2}})}\, }{2\,\sqrt{\rho }\,\, \Re (\zeta^{\prime} ({\displaystyle \frac {1}{2}}  + i\,\rho ))}} \left ( - 
\mathrm{sin}( - \rho \,\mathrm{ln}({\displaystyle \frac {4\,\pi 
}{\sqrt{1 + 4\,\rho ^{2}}}} ) - {\displaystyle \frac {\pi \,
\sqrt{\rho }}{2}}  - \rho  + {\displaystyle \frac {1}{2}} \,
\mathrm{arctan}(2\,\rho ))\right. \\ \left.
\mbox{} + \mathrm{cos}( - \rho \,\mathrm{ln}({\displaystyle 
\frac {4\,\pi }{\sqrt{1 + 4\,\rho ^{2}}}} ) - {\displaystyle 
\frac {\pi \,\sqrt{\rho }}{2}}  - \rho  + {\displaystyle \frac {1
}{2}} \,\mathrm{arctan}(2\,\rho )) - \mathrm{sin}({\displaystyle 
\frac {\pi \,\sqrt{\rho }}{2}} ) + \mathrm{cos}({\displaystyle 
\frac {\pi \,\sqrt{\rho }}{2}} )\right )  {\vrule 
height0.60em width0em depth0.60em}  
 }
}
\label{Q9Asy}
\end{equation}
which places a not-very-strong constraint on the asymptotic limit of $\Re(\zeta^{\prime}(1/2+i\rho))$ as $\rho\rightarrow\infty$, since the sum \eqref{Srule1} must converge (under the weaker requirement that the asymptotic limit of $\zeta(2\rho^2) < \exp(\pi\,\rho)/\rho$ implied in the derivation of \eqref{Srule1}).  
\newline

\subsection{Ratios of $\zeta$}
Since the functional equation and other well-known results in the literature define relationships between ratios of $\zeta$ functions of different argument, it is interesting to explore \eqref{Master1} in this case. In \eqref{Master1} set
\begin{equation}
F(z)= \zeta(4bv(v+i))/\zeta(4av(v+i))
\label{Zrats}
\end{equation}
and for simplicity, let $w=1$ and $(a,b)>0, \in \Re$. Following the approach discussed in the previous section, after translating the contour integral to infinity where it vanishes for $b<a$, and evaluating the residues of the poles so-enclosed we find
\begin{equation}
\mapleinline{inert}{2d}{Sum(Zeta(-b*(2*k+3)*(2*k+1))/Zeta(-a*(2*k+3)*(2*k+1))*(-1)^k,k = 0 ..
infinity)-1/4*Pi*Sum(Zeta(-2*b*n/a)/Zeta(1,-2*n)/(-a^2-2*a*n)^(1/2)/si
nh(1/2*Pi*(-a^2-2*a*n)^(1/2)/a),n = 1 .. infinity) =
(1/4*1/((-b^2+b)^(1/2)*Zeta(a/b)*sinh(1/2*Pi*(-b^2+b)^(1/2)/b))+1/4*Su
m(Zeta(b*tau/a)/Zeta(1,tau)/(-a^2+a*tau)^(1/2)/sinh(1/2*Pi*(-a^2+a*tau
)^(1/2)/a),tau))*Pi+1/2*Zeta(b)/Zeta(a);}{%
\maplemultiline{
  {\displaystyle \sum _{k=0}^{\infty }} \,
{\displaystyle \frac {\zeta ( - b\,K)\,(-1)
^{k}}{\zeta ( - a\,K)}}   
\mbox{} - {\displaystyle \frac {\pi}{4}}
{\displaystyle \sum _{n=1}^{\infty }} \,{\displaystyle \frac {
\zeta ( - {\displaystyle {2\,b\,n}/{a}} )}{\zeta^{\prime} ( - 2
\,n)\,\sqrt{ - a^{2} - 2\,a\,n}\,\mathrm{sinh}({\displaystyle 
\frac {\pi \,\sqrt{ - a^{2} - 2\,a\,n}}{2\,a}} )}} \\  =
  {\vrule height1.71em width0em depth1.71em} 
{\displaystyle \frac {\pi}{4}} \;{\displaystyle \frac {1}{\sqrt{ - 
b^{2} + b}\,\zeta ({\displaystyle \frac {a}{b}} )\,\mathrm{sinh}(
{\displaystyle \frac {\pi \,\sqrt{ - b^{2} + b}}{2\,b}} )}}  
\mbox{} + {\displaystyle \frac {\pi}{4}}\; 
{\displaystyle \sum _{\tau }} \,{\displaystyle \frac {\zeta (
{\displaystyle  {b\,\tau }/{a}} )}{\zeta^{\prime} (\tau )\,\sqrt{
 - a^{2} + a\,\tau }\,\mathrm{sinh}({\displaystyle \frac {\pi \,
\sqrt{ - a^{2} + a\,\tau }}{2\,a}} )}} 
 {\vrule height1.71em width0em depth1.71em} 
\mbox{} + {\displaystyle \frac {1}{2}} \,{\displaystyle \frac {
\zeta (b)}{\zeta (a)}}  }
}
\label{MTx1}
\end{equation}
where
\begin{equation}
K=(2\,k + 3)\,(2\,k + 1)=4(k+1)^2-1
\label{K}
\end{equation}
and both sums on the left-hand side converge only if $b<a$. In the case that $b=2\,m$, a positive even integer, the first sum on the left-hand side of\
\eqref{MTx1} vanishes except for the special case $a=2\,j$ where $j=m+1,\dots$. In that case, a limiting process occurs in both sums on the left-hand side of \eqref{MTx1}, the first occuring between the numerator and denominator of the $k$ sum, the second occuring when the summation index $n$ equals the special value $n=jK$, with K defined by \eqref{K}. The result is
\begin{equation}
\mapleinline{inert}{2d}{-4*Sum(Zeta(1,-2*m*K)*(-1)^k/Zeta(1,-2*m*j*K),k = 0 ..
infinity)/Pi-1/2*Sum(Zeta(-2*m*n/j)/(1+n/j)^(1/2)/Zeta(1,-2*n)/sin(1/2
*Pi*(1+n/j)^(1/2)),n = 1 ..
infinity)-4*m*Sum((-1)^k*Zeta(1,-2*m*K)/Zeta(1,-2*K*j),k = 0 ..
infinity)/Pi =
-2*j/Pi/Zeta(2*j)*Zeta(2*m)-Sum(1/2*Zeta(m*tau/j)/Zeta(1,tau)/sqrt(-1+
1/2*tau/j)/sinh(1/2*Pi*sqrt(-1+1/2*tau/j)),tau)+1/2*j*2^(1/2)/m^(1/2)/
(2*m-1)^(1/2)/sin(1/2*Pi*sqrt(1-1/(2*m)))/Zeta(j/m);}{%
\maplemultiline{
 - {\displaystyle {\frac{4}{\pi}\, {\displaystyle \sum _{k=0
}^{\infty }} \,{\displaystyle \frac {\zeta^{\prime} ( - 2\,m\,K)\,(-1
)^{k}}{\zeta^{\prime} ( - 2\,m\,j\,K)}}  }{}}  - 
{\displaystyle \frac {1}{2}} \, {\displaystyle \sum _{
n=1}^{\infty }}^{\prime} \,{\displaystyle \frac {\zeta ( - {\displaystyle 
 {2\,m\,n}/{j}} )}{\sqrt{1 + {\displaystyle  {n}/{j}} }\,\,
\zeta^{\prime} ( - 2\,n)\,\mathrm{sin} {\displaystyle 
\big{(}\,{\frac{\pi}{2} \,\sqrt{1 + {\displaystyle  {n}/{j}} }}}\;  \big{)}}}    \\
\mbox{} - {\displaystyle { \frac{4\,m}{\pi}\, {\displaystyle 
\sum _{k=0}^{\infty }} \,{\displaystyle \frac {(-1)^{k}\,\zeta^{\prime} (- 2\,m\,K)}{\zeta^{\prime} (- 2\,K\,j)}}  }} 
= 
\mbox{} - {\displaystyle \frac {1}{2}}\, {\displaystyle \sum _{\tau }} 
 \!  \,{\displaystyle \frac {\zeta (
{\displaystyle  {m\,\tau }/{j}} )}{\zeta^{\prime} (\tau )\,\sqrt{
 - 1 + {\displaystyle  {\tau }/{(2\,j)}} }\,\mathrm{sinh}\big{(}\,
{\displaystyle \frac {\pi}{2}}  \,\sqrt{ - 1 + {\displaystyle 
 {\tau }/({2\,j})} }\,\big{)}}}    \\
\mbox{} + {\displaystyle \frac {1}{2}} \,{\displaystyle \frac {j
\,\sqrt{2}}{\sqrt{m}\,\sqrt{2\,m - 1}\,\mathrm{sin}(
{\displaystyle \frac {1}{2}} \,\pi \,\sqrt{1 - {\displaystyle 
 {1}/({2\,m})} })\,\zeta ({\displaystyle  {j}/{m}} )}} - {\displaystyle \frac {2\,j\,\zeta (2\,m)}{\pi \,\zeta (2\,j)}
}   }
}
\label{MtxbrJ}
\end{equation} 
where $\sum^{\prime}$ signifies that all terms satisfying $1+n/j=4h^2$, with $h=1,\dots$ are omitted. The limiting form of these terms is specifically represented by the third sum on the left of \eqref{MtxbrJ}. We note that from a numerical point of view, the two sums indexed by $k$ on the left of \eqref{MtxbrJ} contribute only a small amount to the sum, and that the second sum on the left of this result contains scattered zero elements when $mod(n,j)=0$. Further, the $k>0$ elements of both sums indexed by $k$ on the left of \eqref{MtxbrJ} are, for various choices of $j$ and $m$, at least 25 orders of magnitude smaller than the $k=0$ terms. This, to an excellent degree of approximation, allows us  to re-write the left-hand side of \eqref{MtxbrJ} in the simplified form

\begin{equation}
\mapleinline{inert}{2d}{-4*(-1)^m*((2*Pi)^(6*m*j)*(-1)^(m*j)/Zeta(6*m*j+1)/GAMMA(6*m*j+1)+m*(
-1)^j*(2*Pi)^(6*j)/Zeta(6*j+1)/GAMMA(6*j+1))*Zeta(6*m+1)*GAMMA(6*m+1)/
((2*Pi)^(6*m))/Pi;}{%
\maplemultiline{
 - {\displaystyle \frac {4\,(-1)^{m}\,\left( {\displaystyle \frac {(2\,
\pi )^{(6\,m\,j)}\,(-1)^{(m\,j)}}{\zeta (6\,m\,j + 1)\,\Gamma (6
\,m\,j + 1)}}  + {\displaystyle \frac {m\,(-1)^{j}\,(2\,\pi )^{(6
\,j)}}{\zeta (6\,j + 1)\,\Gamma (6\,j + 1)}} \right) \,\zeta (6\,m + 1)
\,\Gamma (6\,m + 1)}{\pi\,(2\,\pi )^{(6\,m)} }}  \\ -
{\displaystyle \frac {1}{2}} \, {\displaystyle \sum _{
n=1}^{\infty }}^{\prime} \,{\displaystyle \frac {\zeta ( - {\displaystyle 
 {2\,m\,n}/{j}} )}{\sqrt{1 + {\displaystyle  {n}/{j}} }\,\,
\zeta^{\prime} ( - 2\,n)\,\mathrm{sin} {\displaystyle 
\big{(}\,{\frac{\pi}{2} \,\sqrt{1 + {\displaystyle  {n}/{j}} }}}\;  \big{)}}}
}
}
\label{MtxbrJ_a}
\end{equation}
\newline

Any attempt to evaluate the limit $b \rightarrow a$ in \eqref{MTx1} leads to a complicated difference limit of diverging sums on the left-hand side (with a known result obtainable by setting $F(z)=1$ in \eqref{Master1}). However, by comparing the coefficients of the first order expansion of each of the terms in this limit, a new, interesting sum emerges:
\begin{equation}
\mapleinline{inert}{2d}{MtX1 :=
Sum(Zeta(1,-a*(2*k+3)*(2*k+1))*(2*k+3)*(2*k+1)/Zeta(-a*(2*k+3)*(2*k+1)
)*(-1)^k,k = 0 ..
infinity)-1/2*Pi/a*Sum(n/(-a^2-2*a*n)^(1/2)/sinh(1/2*Pi*(-a^2-2*a*n)^(
1/2)/a),n = 1 .. infinity) =
-1/2*Zeta(1,a)/Zeta(a)+1/4*Pi/(-a*(-1+a))^(1/2)/a/sinh(1/2*Pi*(-a^2+a)
^(1/2)/a)-1/4*Pi/a*Sum(tau/(-a^2+a*tau)^(1/2)/sinh(1/2*Pi*(-a^2+a*tau)
^(1/2)/a),tau);}{%
\maplemultiline{
{\displaystyle \sum _{k=0}^{\infty }
} \,{\displaystyle \frac {K\,\zeta^{\prime} (- a\,K)\,(-1)^{k}}{\zeta ( - a\,K)}}    
\mbox{} - {\displaystyle \frac {\pi}{2\,a}} \,{\displaystyle {
{\displaystyle \sum _{n=1}^{\infty }} \,
{\displaystyle \frac {n}{\sqrt{ - a^{2} - 2\,a\,n}\,\mathrm{sinh}
({\displaystyle \frac {\pi \,\sqrt{ - a^{2} - 2\,a\,n}}{2\,a}} )}
}  }{}}\\ = - {\displaystyle \frac {1}{2}} \,
{\displaystyle \frac {\zeta^{\prime} (a)}{\zeta (a)}}  
\mbox{} + {\displaystyle \frac {1}{4}} \,{\displaystyle \frac {
\pi }{\sqrt{ - a\,( - 1 + a)}\,a\,\mathrm{sinh}({\displaystyle 
\frac {\pi \,\sqrt{ - a^{2} + a}}{2\,a}} )}}  
\mbox{} - {\displaystyle \frac {\pi}{4\,a}} \,{\displaystyle {
{\displaystyle \sum _{\tau }} \,{\displaystyle 
\frac {\tau }{\sqrt{ - a^{2} + a\,\tau }\,\mathrm{sinh}(
{\displaystyle \frac {\pi \,\sqrt{ - a^{2} + a\,\tau }}{2\,a}} )}
}   }{}}  }
}
\label{Mtx1a=b}
\end{equation}

valid for $a<0$. From the properties of the $\zeta$ function, it is possible to show that for large values of $k$ with $a<0$, the terms of the first sum on the left-hand side of \eqref{Mtx1a=b} asymptotically are approximated by
\begin{equation}
\mapleinline{inert}{2d}{Asy :=
Zeta(1,-a*(2*k+3)*(2*k+1))*(2*k+3)*(2*k+1)/Zeta(-a*(2*k+3)*(2*k+1)) =
-(2*k+3)*(2*k+1)*ln(2)/(2^(-a*(2*k+3)*(2*k+1)));}{%
\[
{\displaystyle \frac {K\,\zeta^{\prime} ( - a\,K)}{\zeta ( - a\,K))}} \approx - {\displaystyle \frac {
\mathrm{ln}(2^{K})}{2^{( - a\,K)}}} 
\]
\label{Kasy}
}
\end{equation}
so the corresponding series in \eqref{Mtx1a=b} converges.  Notice that neither the term involving the sum over zeros of $\zeta(\tau)$ nor the sum over the index $n$ in \eqref{Mtx1a=b} involve $\zeta(\tau)$ or any of its derivatives, in contrast to previous results. In the case that $a=-1/K$ where $K$ labels one index of the first sum in \eqref{Mtx1a=b}, a limit exists between the divergent term and the second term on the right of that equation. For example, in the case $a=-1/3$, we find

\begin{equation}
\mapleinline{inert}{2d}{Sum(tau/(-3*tau-1)^(1/2)/sinh(1/2*Pi*(-3*tau-1)^(1/2)),tau) =
-4/9*Sum(K*Zeta(1,1/3*K)*(-1)^k/Zeta(1/3*K),k = 1 ..
infinity)/Pi+2*Sum(n/(6*n-1)^(1/2)/sinh(1/2*Pi*(6*n-1)^(1/2)),n = 1 ..
infinity)+4/9*(39/16-3*gamma)/Pi-2/9*Zeta(1,-1/3)/Pi/Zeta(-1/3);}{%
\maplemultiline{
{\displaystyle \sum _{\tau }} \,{\displaystyle \frac {\tau }{
\sqrt{ - 3\,\tau  - 1}\,\mathrm{sinh}({\displaystyle \frac {\pi 
\,\sqrt{ - 3\,\tau  - 1}}{2}} )}} = - {\displaystyle \frac {4}{9\,\pi}
} \,{\displaystyle {{\displaystyle \sum _{k=1}^{\infty }} 
\,{\displaystyle \frac {K\,\zeta^{\prime} ({\displaystyle  {K}/{3
}} )\,(-1)^{k}}{\zeta ({\displaystyle {K}/{3}} )}} }{}} 
 \\
\mbox{} + 2 {\displaystyle \sum _{n=1}^{\infty }} \,
{\displaystyle \frac {n}{\sqrt{6\,n - 1}\,\mathrm{sinh}(
{\displaystyle \frac {\pi \,\sqrt{6\,n - 1}}{2}} )}}   + {\displaystyle \frac{4}{9\,\pi}\,\displaystyle  {({\displaystyle \frac {39}{
16}}  - 3\,\gamma )}{ }}  - {\displaystyle \frac {2}{9}} \,
{\displaystyle \frac {\zeta^{\prime} ({\displaystyle {-1}/{3}} )
}{\pi \,\zeta ({\displaystyle {-1}/{3}} )}}  \,.}
}
\label{eq4p25a}
\end{equation}

As well, the first sum in \eqref{Mtx1a=b} can be written in the interesting form
\begin{equation}
  \mapleinline{inert}{2d}{Diff(ln(product(Zeta(-a*(4*k+1)*(4*k+3))/Zeta(-a*(4*k+3)*(4*k+5)),k =
0 .. infinity)),a);}{%
\[ {\displaystyle \sum _{k=0}^{\infty }
} \,{\displaystyle \frac {K\,\zeta^{\prime} (- a\,K)\,(-1)^{k}}{\zeta ( - a\,K)}}\,=\,
{\displaystyle \frac {\partial }{\partial a}}\,\mathrm{ln} \left(  \! 
{\displaystyle \prod _{k=0}^{\infty }} \,{\displaystyle \frac {
\zeta ( - a\,(4\,k + 1)\,(4\,k + 3))}{\zeta ( - a\,(4\,k + 3)\,(4
\,k + 5))}}  \!  \right) 
\]
}
\label{diffLog}
  \end{equation}

\subsection{Derivative Families of Sum Rules}

From \eqref{Master1_Zeta}, new families of sum rules can be obtained. After the change of variables $v=i\,t$ in \eqref{Master1_Zeta}, translate the resulting integration contour to the right by a small amount $\delta$, and, provided that the integrand contains no singularities for $0\leq\delta$, the value of the integral does not change and we have 
\begin{equation}
\mapleinline{inert}{2d}{MT2b :=
Int(w^(-4*a*(t+delta)*(t+delta-1))/Zeta(-4*a*(t+delta)*(t+delta-1))/co
s(Pi*(t+delta)),t = -i*infinity .. i*infinity) = i*w^a/Zeta(a);}{%
\[
{\displaystyle \int _{ - i\,\infty }^{i\,\infty 
}} {\displaystyle \frac {w^{( - 4\,a\,(t + \delta )\,(t + \delta 
 - 1))}}{\zeta ( - 4\,a\,(t + \delta )\,(t + \delta  - 1))\,
\mathrm{cos}(\pi \,(t + \delta ))}} \,dt={\displaystyle \frac {i
\,w^{a}}{\zeta (a)}} 
\]
}
\label{MT2b}
\end{equation}
or, equivalently
\begin{equation}
\mapleinline{inert}{2d}{Diff(f(delta),delta);}{%
\[
{\displaystyle \frac {d}{d\delta }}\,{\displaystyle \int _{ - i\,\infty }^{i\,\infty 
}} {\displaystyle \frac {w^{( - 4\,a\,(t + \delta )\,(t + \delta 
 - 1))}}{\zeta ( - 4\,a\,(t + \delta )\,(t + \delta  - 1))\,
\mathrm{cos}(\pi \,(t + \delta ))}} \,dt=0\,.
\]
}
\label{diffdelta}
\end{equation}

To zero'th order in $\delta$ as $\delta\rightarrow0$ after setting $w=1$, (and equivalent to integration by parts), we find, for $\Re(a)>0$, the identity
\begin{equation}
\mapleinline{inert}{2d}{MT3a := Int(sinh(Pi*v)/Zeta(-4*i*a*v*(i*v-1))/cosh(Pi*v)^2,v =
-infinity .. infinity) =
4*i*a/Pi*Int((2*i*v-1)*Zeta^{\prime}(-4*i*a*v*(i*v-1))/cosh(Pi*v)/Zeta(-4*i*a
*v*(i*v-1))^2,v = -infinity .. infinity);}{%
\maplemultiline{
{\displaystyle \int _{ - \infty }^{\infty }} 
{\displaystyle \frac {\mathrm{sinh}(\pi \,v)}{\zeta ( - 4\,i\,a\,
v\,(i\,v - 1))\,\mathrm{cosh}(\pi \,v)^{2}}} \,dv= 
{\displaystyle \frac {4\,i\,a}{\pi }} \,{\displaystyle \int _{ - 
\infty }^{\infty }} {\displaystyle \frac {(2\,i\,v - 1)\,\zeta^{\prime} ( - 4\,i\,a\,v\,(i\,v - 1))}{\mathrm{cosh}(\pi \,v)\,\zeta (
 - 4\,i\,a\,v\,(i\,v - 1))^{2}}} \,dv } \,.
}
\label{MT3a}
\end{equation}
Unfortunately, the kernels in both integrals in \eqref{MT3a} do not satisfy the requirements for a new master equation of the type \eqref{Master1}, (see \cite[Eq. (2.5)] {Master} ), so neither can be evaluated in analogy to \eqref{Master1}. However, by translating the integration contour of both integrals parallel to the real $v$ axis to $\Im(v) =\infty$ as in the previous sections, and evaluating the residues of the poles encountered along the way, we find a new relationship among the zeros of the $\zeta$ function:
\begin{equation}
\mapleinline{inert}{2d}{Sum((1/2*coth(1/2*Pi*(-a+tau)^(1/2)/a^(1/2))+a^(1/2)/Pi/(-a+tau)^(1/2
))/Zeta(1,tau)/sinh(1/2*Pi*(-a+tau)^(1/2)/a^(1/2))/(a-tau),tau) =
4*a*Sum((-1)^k*Zeta(1,-a*(2*k+3)*(2*k+1))/Zeta(-a*(2*k+3)*(2*k+1))^2,k
= 0 ..
infinity)/Pi^2+Sum((1/2*a^(1/2)/Pi/(a+2*n)^(1/2)+1/4*cot(1/2*Pi*(a+2*n
)^(1/2)/a^(1/2)))/Zeta(1,-2*n)/sin(1/2*Pi*(a+2*n)^(1/2)/a^(1/2))/(a+2*
n),n = 1 .. infinity)-2*a*Zeta(1,a)/Pi^2/Zeta(a)^2;}{%
\maplemultiline{
{\displaystyle \sum _{\tau }} \,{\displaystyle \frac {
{\displaystyle \frac {1}{2}} \,\mathrm{coth}({\displaystyle 
\frac {\pi \,\sqrt{ - a + \tau }}{2\,\sqrt{a}}} ) + 
{\displaystyle \frac {\sqrt{a}}{\pi \,\sqrt{ - a + \tau }}} }{
\zeta^{\prime} (\tau )\,\mathrm{sinh}({\displaystyle \frac {\pi \,
\sqrt{ - a + \tau }}{2\,\sqrt{a}}} )\,(a - \tau )}} =
{\displaystyle {\displaystyle \frac{4\,a}{\pi^{2}}}  {\displaystyle \sum _{k=0
}^{\infty }} \,{\displaystyle \frac {(-1)^{k}\,\zeta^{\prime} ( - a\,
K)}{\zeta ( - a\,K)^{2}
}}  }  - 
{\displaystyle \frac {2\,a\,\zeta^{\prime} (a)}{\pi ^{2}\,\zeta (a)^{
2}}}  \\
\hspace{6.5cm} +  {\displaystyle \sum _{n=1}^{\infty }} \,
{\displaystyle \frac {{\displaystyle \frac {\sqrt{a}}{2\,\pi \,
\sqrt{a + 2\,n}}}  + {\displaystyle \frac {1}{4}} \,\mathrm{cot}(
{\displaystyle \frac {\pi \,\sqrt{a + 2\,n}}{2\,\sqrt{a}}} )}{
\zeta^{\prime} ( - 2\,n)\,\mathrm{sin}({\displaystyle \frac {\pi \,
\sqrt{a + 2\,n}}{2\,\sqrt{a}}} )\,(a + 2\,n)}}  \,. } 
}
\label{S4}
\end{equation}    

\subsection{Application of the Functional Equation}

Following the motivation of \eqref{Zrats} for $a>0$, we choose
\begin{equation}
F(z)= \zeta(4av(v+i))/\zeta(1-4av(v+i))
\label{ZratFunc}
\end{equation}
and, making allowance for singularities existing in $\mathfrak{S}$ we obtain 
\begin{equation}
\mapleinline{inert}{2d}{Int(w^(v*(v+i))*Zeta(4*a*v*(v+i))/Zeta(1-4*a*v*(v+i))/cosh(Pi*v),v =
-infinity .. infinity) =
w^(1/4)*Zeta(a)/Zeta(1-a)-Sum(1/2*Pi*w^(1/2*n/a)*w^(1/(4*a))*Zeta(2*n+
1)/Zeta(1,-2*n)/(-a^2+2*a*n+a)^(1/2)/sinh(1/2*Pi*(-a^2+2*a*n+a)^(1/2)/
a),n = 1 ..
infinity)-Pi*w^(1/(4*a))/(-a^2+a)^(1/2)/sinh(1/2*Pi*(-a^2+a)^(1/2)/a);
}{%
\maplemultiline{
{\displaystyle \int _{ - \infty }^{\infty }} {\displaystyle 
\frac {w^{(v\,(v + i))}\,\zeta (4\,a\,v\,(v + i))}{\zeta (1 - 4\,
a\,v\,(v + i))\,\mathrm{cosh}(\pi \,v)}} \,dv={\displaystyle 
\frac {w^{(1/4)}\,\zeta (a)}{\zeta (1 - a)}} - {\displaystyle \frac {\pi \,w^{( {1}/{4\,a})}}{
\sqrt{ - a^{2} + a}\,\mathrm{sinh}({\displaystyle \frac {\pi \,
\sqrt{ - a^{2} + a}}{2\,a}} )}} \\
\mbox{} -  {\displaystyle \frac {\pi}{2}} {\displaystyle \sum _{n=1}^{N }} \,
 \,{\displaystyle \frac {
{\displaystyle w^{ (\frac {2\,n+1}{4\,a})}}\,\zeta (2\,n
 + 1)}{\zeta^{\prime}( - 2\,n)\,\sqrt{ - a^{2} + 2\,a\,n + a}\,
\mathrm{sinh}({\displaystyle \frac {\pi \,\sqrt{ - a^{2} + 2\,a\,
n + a}}{2\,a}} )}}  \;.  }
}
\label{Zrat1}
\end{equation}

The result \eqref{Zrat1} is valid for $a>0,a \in \mathfrak{R}$ with $N\rightarrow \infty$. In the case that $a\in\mathfrak{C}$ the value of $N$ must be chosen as discussed in previous sections to include only those poles lying in $\mathfrak{S}$ for particular values of $a$. Notice that all residues attached to zeros of the term $\zeta (1 - 4\,
a\,v\,(v + i))$ in the denominator of \eqref{Zrat1} are respectively proportional to $\zeta(1-\tau)$ and therefore vanish.  Applying the functional equation for zeta functions to \eqref{Zrat1} and employing well-known duplication and inversion formula for $\Gamma$ functions, eventually yields

\begin{equation}
\mapleinline{inert}{2d}{Int(w^(v*(v+i))*(2*Pi)^(4*a*v*(v+i))/cos(2*Pi*a*v*(v+i))/GAMMA(4*a*v*
(v+i))/cosh(Pi*v),v = -infinity .. infinity) =
w^(1/4)*(2*Pi)^a/cos(1/2*Pi*a)/GAMMA(a)-2*Pi*Sum((-1)^n*w^(1/4*(2*n+1)
/a)*(2*Pi)^(2*n)/GAMMA(2*n+1)/(-a^2+2*a*n+a)^(1/2)/sinh(1/2*Pi*(-a^2+2
*a*n+a)^(1/2)/a),n = 0 .. N);}{%
\maplemultiline{
{\displaystyle \int _{ - \infty }^{\infty }} {\displaystyle 
\frac {w^{(v\,(v + i))}\,(2\,\pi )^{(4\,a\,v\,(v + i))}}{\mathrm{
cos}(2\,\pi \,a\,v\,(v + i))\,\Gamma (4\,a\,v\,(v + i))\,\mathrm{
cosh}(\pi \,v)}} \,dv= 
{\displaystyle \frac {w^{(1/4)}\,(2\,\pi )^{a}}{\mathrm{cos}(
{\displaystyle \frac {\pi \,a}{2}} )\,\Gamma (a)}} \\ - 2\,\pi \,
  {\displaystyle \sum _{n=0}^{N}} \,{\displaystyle 
\frac {(-1)^{n}\,w^{(\frac {2\,n + 1}{4\,a})}\,(2\,\pi )^{(2\,n)}
}{\Gamma (2\,n + 1)\,\sqrt{ - a^{2} + 2\,a\,n + a}\,\mathrm{sinh}
({\displaystyle \frac {\pi \,\sqrt{ - a^{2} + 2\,a\,n + a}}{2\,a}
} )}} \;, }
}
\label{Zrtag3}
\end{equation}
from which all reference to $\zeta(\tau)$ has vanished. We note that the integrand of \eqref{Zrtag3} obeys \eqref{f(v)}, and therefore exists as a master kernel in its own right, in analogy to \eqref{Master1}, from which \eqref{Zrtag3} could have been otherwise obtained. In contrast to previous sections, the integrands of both integrals in \eqref{Zrat1} and \eqref{Zrtag3} diverge as the integration contour is translated to $i \,\infty$, which precludes the possibility of replacing the integrals by a sum over the residues that may have been encountered. 
 
\subsection{Application of alternate Master equations}

As stated in section 1, any integral with infinite limits, whose integrand obeys \eqref{f(v)} and has appropriate asymptotic properties, can be analytically evaluated using the analogue of \eqref{Master1}, provided the residues of poles residing in $\mathfrak{S}$ are added to the right-hand side. In this example, we shall use
\begin{equation}
\mapleinline{inert}{2d}{F(z) = 1/(Zeta(a*(v+i))*Zeta(-a*v));}{%
\[
\mathrm{F}(z)={\displaystyle \frac {1}{\zeta (a\,(v + i))\,\zeta 
( - a\,v)}} 
\]
}
\label{Zeta^2}
\end{equation}

in \eqref{1/Master1}, set $w=1$, assume $a=\alpha+i\beta$ where necessary, and, with the requirement that $\beta\neq 0$ we obtain
\begin{equation}
\mapleinline{inert}{2d}{Int(1/(Zeta(a*(v+i))*Zeta(-a*v)*cosh(Pi*v)),v = -infinity ..
infinity) =
1/(Zeta(1/2*i*a)^2)-2*i*Pi/a*Sum(-1/(Zeta(1,-2*n)*Zeta(i*a+2*n)*cosh(2
*Pi*n/a)),n = 1 ..
N)-2*i*Pi/a*Sum(-1/(Zeta(1,tau)*Zeta(i*a-tau)*cosh(Pi*tau/a)),`Im(tau)
>0`);}{%
\maplemultiline{
{\displaystyle \int _{ - \infty }^{\infty }} {\displaystyle 
\frac {1}{\zeta (a\,(v + i))\,\zeta ( - a\,v)\,\mathrm{cosh}(\pi 
\,v)}} \,dv={\displaystyle \frac {1}{\zeta ({\displaystyle 
\frac {i\,a}{2}} )^{\displaystyle ^{2}}}}  
\mbox{} + {\displaystyle  {\frac{2\,i\,\pi}{a} 
{\displaystyle \sum _{n=1}^{N}}  {\displaystyle 
\frac {1}{\zeta^{\prime} ( - 2\,n)\,\zeta (i\,a + 2\,n)\,\mathrm{cosh
}({\displaystyle \frac {2\,\pi \,n}{a}} )}} 
 }}  \\
\mbox{} + {\displaystyle {\frac{2\,i\,\pi}{a} 
{\displaystyle \sum _{\mathit{\Im(\tau)>0}}}   
{\displaystyle \frac {1}{\zeta^{\prime} (\tau )\,\zeta (i\,a - \tau )
\,\mathrm{cosh}({\displaystyle {\pi \,\tau }/{a}} )}} 
 }{}}  \, .}
}
\label{M1z^2}
\end{equation}

The terms appearing on the right-hand side of \eqref{M1z^2} follow in exact analogy to similar terms appearing in previous examples. The first term corresponds to the residue of $1/\cosh(\pi v$) as it appears in \eqref{Master1}, the second sum is limited by the upper index $N$ to any residues belonging to the trivial zeros $\zeta(-2n)$ that may reside in $\mathfrak{S}$ in the case that $\beta>0$ and the last term corresponds to similar residues that may need to be included in the case $\beta<0$. As an example, in the case $a=1-4i$ we find
\begin{equation}
\mapleinline{inert}{2d}{Int(1/(Zeta((1-4*i)*(v+i))*Zeta((-1+4*i)*v)*cosh(Pi*v)),v = -infinity
.. infinity) =
1/(Zeta(2+1/2*i)^2)+(-8/17+2/17*i)*Pi/Zeta(1,1/2+i*rho[1])/Zeta(7/2+i-
i*rho[1])/cosh((1/34+2/17*i)*Pi*(1+2*i*rho[1]));}{%
\maplemultiline{
{\displaystyle \int _{ - \infty }^{\infty }} {\displaystyle 
\frac {1}{\zeta ((1 - 4\,i)\,(v + i))\,\zeta (( - 1 + 4\,i)\,v)\,
\mathrm{cosh}(\pi \,v)}} \,dv= 
{\displaystyle \frac {1}{\zeta (2 + {\displaystyle {i}/{2}} 
)^{\displaystyle ^{2}}}} \\ +\; {\displaystyle \frac {2\,\pi\,( - {\displaystyle {4}{}
}  + {\displaystyle  {i}{}} )/17 }{\zeta^{\prime} (
{\displaystyle  {1}/{2}}  + i\,{\rho _{1}})\,\zeta (
{\displaystyle  {7}/{2}}  + i - i\,{\rho _{1}})\,\mathrm{cosh
}(({\displaystyle  {1}/{34}}  + {\displaystyle {2\,i}/{
17}} )\,\pi \,(1 + 2\,i\,{\rho _{1}}))}}  \, .}
}
\label{ex4p6}
\end{equation}
Again, following the procedure used in previous sections, we translate the integration contour to $v=i\infty$, adding residues of poles encountered, and after some re-arrangement, eventually arrive at the sum rule
\begin{equation}
\mapleinline{inert}{2d}{Sum(1/(Zeta(1,tau)*a*Zeta(i*a-tau)*cosh(Pi*tau/a)),tau) =
-1/a*Sum(1/(Zeta(1,-2*n)*Zeta(i*a+2*n)*cosh(2*Pi*n/a)),n = 1 ..
infinity)-i/Pi*Sum((-1)^k/Zeta(i*a*k+3/2*i*a)/Zeta(-i*a*k-1/2*i*a),k =
0 .. infinity)+1/2*i/Pi/Zeta(1/2*i*a)^2;}{%
\maplemultiline{
{\displaystyle \sum _{\tau }} \,{\displaystyle \frac {1}{\zeta^{\prime} (\tau )\,\zeta (i\,a - \tau )\,\mathrm{cosh}(
{\displaystyle  {\pi \,\tau }/{a}} )}} = - {\displaystyle 
 {{\displaystyle \sum _{n=1}^{\infty }} \,{\displaystyle 
\frac {1}{\zeta^{\prime} ( - 2\,n)\,\zeta (i\,a + 2\,n)\,\mathrm{cosh
}({\displaystyle \frac {2\,\pi \,n}{a}} )}} }{}}  \\
\mbox{} - {\displaystyle  {\displaystyle \frac {i\,a}{\pi}\,{\displaystyle 
\sum _{k=0}^{\infty }} \,{\displaystyle \frac {(-1)^{k}}{\zeta (i
\,a\,k + {\displaystyle  {3}} \,i\,a/2)\,\zeta ( - i\,a\,k
 -  \,i\,a/2)}} }{ }} 
 + \,{\displaystyle \frac {i\,a}{2\,\pi \,
\zeta ({\displaystyle {i\,a}/{2}} )^{2}}}\; .  }
}
\label{Mt1a}
\end{equation}

It is worth noting that for small values of $\alpha$ and $\beta$, the three component terms of the right-hand side of \eqref{Mt1a} undergo severe numerical cancellation of significant digits, and the left-hand side typically needs only one term to yield numerical equality to several digits. As examples, we have, after some manipulation involving the functional equation for the $\zeta$ function, duplication and reflection properties of the $\Gamma$-function and the invocation of \eqref{ZetaprimeNegN} the following results. For $a=-i$ it turns out that residues belonging to complex conjugate values of the summation index $\tau$ cancel one another, so that the left-hand side of \eqref{Mt1a} vanishes, leading to the following transformation rule between sums
\begin{equation}
\mapleinline{inert}{2d}{2^(1/2)*Pi^(1/2)*Sum((-2*Pi)^k/Zeta(k+3/2)^2/sin(1/2*Pi*(k+1/2))/GAMM
A(k+3/2),k = 1 .. infinity) =
-4/Zeta(3/2)^2-1/2*1/(Zeta(1/2)^2*Pi)+2*Sum((-1)^n*(2*Pi)^(2*n)/Zeta(1
+2*n)^2/GAMMA(1+2*n),n = 1 .. infinity);}{%
\maplemultiline{
\sqrt{2\,\pi }\, {\displaystyle \sum _{k=1}^{
\infty }} \,{\displaystyle \frac {( - 2\,\pi )^{k}}{\zeta (k + 
{\displaystyle  {3}/{2}} )^{\displaystyle^{2}}\,\mathrm{sin}\, \left(  
{\displaystyle  {(2\,k + {\displaystyle 1} )\,\pi/4}{}}    \right) \,\Gamma (k + {\displaystyle  {3}/{2}} )}} 
 = 
 - {\displaystyle \frac {4}{\zeta ({\displaystyle  {3}/{2}} )
^{2}}}  - {\displaystyle \frac {1}{2\,\pi}} \,{\displaystyle \frac {1
}{\zeta ({\displaystyle  {1}/{2}} )^{2}}} \\\hspace{8.3cm} + 2\,
  {\displaystyle \sum _{n=1}^{\infty }} \,
{\displaystyle \frac {(-1)^{n}\,(2\,\pi )^{(2\,n)}}{\zeta (1 + 2
\,n)^{2}\,\Gamma (1 + 2\,n)}} \,.  }
}
\label{a=-I}
\end{equation}
For the case $a=2\,i$ we find
\begin{equation}
\mapleinline{inert}{2d}{Mta4 := Sum(1/(Zeta(1,tau)*Zeta(-2-tau)*cos(1/2*Pi*tau)),tau) =
-1/(Pi*Zeta(-1)^2)-2/Zeta(1,-2)-8*Pi^2*Sum((2*Pi)^(2*k)*((2*k+3)/Zeta(
2*k+3)/Zeta(2*k)-2*Pi/Zeta(4+2*k)/Zeta(2*k+1))/GAMMA(4+2*k),k = 1 ..
infinity);}{%
\maplemultiline{
 {\displaystyle \sum _{\tau }} \,{\displaystyle 
\frac {1}{\zeta^{\prime} (\tau )\,\zeta ( - 2 - \tau )\,\mathrm{cos}(
{\displaystyle {\pi \,\tau }/{2}} )}} = - {\displaystyle 
\frac {1}{\pi \,\zeta (-1)^{2}}}  - {\displaystyle \frac {2}{
\zeta^{\prime} (-2)}}  \\
\mbox{} - 8\,\pi ^{2}\,{\displaystyle \sum _{k=1}^{
\infty }} \,{\displaystyle \frac {(2\,\pi )^{2\,k}\,\left(
{\displaystyle \frac {2\,k + 3}{\zeta (2\,k + 3)\,\zeta (2\,k)}} 
 - {\displaystyle \frac {2\,\pi }{\zeta (4 + 2\,k)\,\zeta (2\,k
 + 1)}} \right) }{\Gamma (4 + 2\,k)}}  }
}
\label{Mta4}
\end{equation}
and for $a=i$ we obtain
\begin{equation}
\mapleinline{inert}{2d}{Mt5b := Sum(1/(Zeta(1,tau)*Zeta(-1-tau)*cos(Pi*tau)),tau) =
8*Pi^2*Sum((-1)^k*(2*Pi)^(2*k)/Zeta(2*k+3)/GAMMA(2*k+3)/Zeta(2*k+1),k
= 1 ..
infinity)+1/(Zeta(-3/2)*Zeta(1/2)*Pi)+Pi*Sum((-Pi)^k*GAMMA(-1/2*k-3/4)
/GAMMA(5/4+1/2*k)/Zeta(5/2+k)/Zeta(k+1/2),k = 1 ..
infinity)-1/2*1/(Pi*Zeta(-1/2)^2);}{%
\maplemultiline{
{\displaystyle \sum _{\tau }} \,{\displaystyle 
\frac {1}{\zeta^{\prime}(\tau )\,\zeta ( - 1 - \tau )\,\mathrm{cos}(
\pi \,\tau )}} =8\,\pi ^{2}\, {\displaystyle \sum _{k=
1}^{\infty }} \,{\displaystyle \frac {(-1)^{k}\,(2\,\pi )^{(2\,k)
}}{\zeta (2\,k + 3)\,\Gamma (2\,k + 3)\,\zeta (2\,k + 1)}}    \\
\mbox{} + {\displaystyle \frac {1}{\zeta ({\displaystyle  {
-3}/{2}} )\,\zeta ({\displaystyle {1}/{2}} )\,\pi }}  + \pi 
\, {\displaystyle \sum _{k=1}^{\infty }} \,
{\displaystyle \frac {( - \pi )^{k}\,\Gamma ( - {\displaystyle 
 {k}/{2}}  - {\displaystyle  {3}/{4}} )}{\Gamma (
{\displaystyle  {5}/{4}}  + {\displaystyle  {k}/{2}} )\,
\zeta ({\displaystyle  {5}/{2}}  + k)\,\zeta (k + 
{\displaystyle {1}/{2}} )}}    - {\displaystyle 
\frac {1}{2\,\pi}} \,{\displaystyle \frac {1}{\zeta (
{\displaystyle  {-1}/{2}} )^{2}}} \,. }
}
\label{Mt5b}
\end{equation}
Continuing in the same vein, the case $a=-i/2$ yields the interesting result
\begin{equation}
\mapleinline{inert}{2d}{Mt3x := Sum(1/(Zeta(1,tau)*Zeta(1/2-tau)*cos(2*Pi*tau)),tau) =
-2*Sum((-1)^k*(2*Pi)^(2*k)/Zeta(1/2+2*k)/Zeta(2*k+1)/GAMMA(2*k+1),k =
1 .. infinity)-1/2*Sum((-1)^k/Zeta(1/2*k+3/4)/Zeta(-1/2*k-1/4),k = 0
.. infinity)/Pi+1/4*1/(Pi*Zeta(1/4)^2);}{%
\maplemultiline{
{\displaystyle \sum _{\tau }} \,{\displaystyle 
\frac {1}{\zeta^{\prime}(\tau )\,\zeta ({\displaystyle \frac {1}{2}
}  - \tau )\,\mathrm{cos}(2\,\pi \,\tau )}} =
 -  {\displaystyle 2\, \sum _{k=1}^{\infty }} 
 \; {\displaystyle \frac {(-1)^{k}\,(2\,\pi) ^{(2\,k)}}{
\zeta ({\displaystyle  {1}/{2}}  + 2\,k)\,\zeta (2\,k + 1)\,
\Gamma (2\,k + 1)}}   \\ - {\displaystyle 
\frac {1}{2\,\pi}} \,{\displaystyle  {{\displaystyle \sum _{k=0}
^{\infty }} \,{\displaystyle \frac {(-1)^{k}}{\zeta (
{\displaystyle{k}/{2}}  + {\displaystyle {3}/{4}} )\,
\zeta ( - {\displaystyle  {k}/{2}}  - {\displaystyle  {1
}/{4}} )}} }{ }}  
\mbox{} + {\displaystyle \frac {1}{4\,\pi}} \,{\displaystyle \frac {1
}{\zeta ({\displaystyle  {1}/{4}} )^{2}}}  }
}
\label{Mt3x}
\end{equation}
which, to an excellent degree of approximation, can be written
\begin{equation}
\mapleinline{inert}{2d}{(-1/(Zeta(1,1/2+i*rho[1])*Zeta(-i*rho[1]))-1/(Zeta(1,1/2-i*rho[1])*Ze
ta(i*rho[1])))/cosh(2*Pi*rho[1]) =
-Sum(2*(-1)^k*4^k*Pi^(2*k)/Zeta(1/2+2*k)/Zeta(2*k+1)/GAMMA(2*k+1),k =
1 .. infinity)-1/2*Sum((-1)^k/Zeta(1/2*k+3/4)/Zeta(-1/2*k-1/4),k = 0
.. infinity)/Pi+1/4*1/(Pi*Zeta(1/4)^2);}{%
\maplemultiline{
{\displaystyle \frac { - {\displaystyle \frac {1}{\zeta^{\prime}(
{\displaystyle  {1}/{2}}  + i\,{\rho _{1}})\,\zeta ( - i\,{
\rho _{1}})}}  - {\displaystyle \frac {1}{\zeta^{\prime}(
{\displaystyle {1}/{2}}  - i\,{\rho _{1}})\,\zeta (i\,{\rho 
_{1}})}} }{\mathrm{cosh}(2\,\pi \,{\rho _{1}})}} \approx 
{\displaystyle \frac {1}{4\,\pi}} \,{\displaystyle \frac {1
}{\zeta ({\displaystyle  {1}/{4}} )^{2}}} \\
 -  {\displaystyle 2\, \sum _{k=1}^{\infty }} 
 \; {\displaystyle \frac {(-1)^{k}\,(2\,\pi) ^{(2\,k)}}{
\zeta ({\displaystyle  {1}/{2}}  + 2\,k)\,\zeta (2\,k + 1)\,
\Gamma (2\,k + 1)}}    - {\displaystyle 
\frac {1}{2\,\pi}} \,{\displaystyle  {{\displaystyle \sum _{k=0}
^{\infty }} \,{\displaystyle \frac {(-1)^{k}}{\zeta (
{\displaystyle{k}/{2}}  + {\displaystyle {3}/{4}} )\,
\zeta ( - {\displaystyle  {k}/{2}}  - {\displaystyle  {1
}/{4}} )}} }{ }}  
\mbox{}   }
}
\label{Mt3ax}
\end{equation}
because the contribution of the second, and higher values of $\tau$ in the left-hand sum of \eqref{Mt3x} is at least 20 orders of magnitude smaller than the first, due to the reciprocal term   $1/\cosh(2\pi\rho_{1})$ appearing in \eqref{Mt3ax}. For the case $a=i/2$ we find the closely related result
\begin{equation}
\mapleinline{inert}{2d}{(-1/(Zeta(1,1/2+i*rho[1])*Zeta(-1-i*rho[1]))-1/(Zeta(1,1/2-i*rho[1])*
Zeta(-1+i*rho[1])))/cosh(2*Pi*rho[1]) =
-2*Sum((-1)^k*(2*Pi)^(2*k)/Zeta(2*k+1)/GAMMA(2*k+1)/Zeta(-1/2+2*k),k =
1 .. infinity)+1/2*Sum((-1)^k/Zeta(-1/2*k-3/4)/Zeta(1/2*k+1/4),k = 0
.. infinity)/Pi-1/4*1/(Pi*Zeta(-1/4)^2);}{%
\maplemultiline{
{\displaystyle \frac { - {\displaystyle \frac {1}{\zeta^{\prime}(
{\displaystyle  {1}/{2}}  + i\,{\rho _{1}})\,\zeta ( - 1 - i
\,{\rho _{1}})}}  - {\displaystyle \frac {1}{\zeta^{\prime}(
{\displaystyle  {1}/{2}}  - i\,{\rho _{1}})\,\zeta ( - 1 + i
\,{\rho _{1}})}} }{\mathrm{cosh}(2\,\pi \,{\rho _{1}})}} \approx  
- {\displaystyle \frac {1}{4\,\pi}} \,{\displaystyle \frac {1
}{\zeta ({\displaystyle  {-1}/{4}} )^{2}}} \\
 - 2\, {\displaystyle \sum _{k=1}^{\infty }} \,
{\displaystyle \frac {(-1)^{k}\,(2\,\pi )^{(2\,k)}}{\zeta ( - {\displaystyle {1}/{2}} 
 + 2\,k)\,\zeta (2\,k
 + 1)\,\Gamma (2\,k + 1)\,}}  + {\displaystyle \frac {1}{2\,\pi}} \,
{\displaystyle {{\displaystyle \sum _{k=0}^{\infty }} \,
{\displaystyle \frac {(-1)^{k}}{\zeta ( - {\displaystyle {k
}/{2}}  - {\displaystyle  {3}/{4}} )\,\zeta ({\displaystyle 
 {k}/{2}}  + {\displaystyle {1}/{4}} )}} }{ }}  
\mbox{} }
}
\label{Mt4}
\end{equation}

Unfortunately, the right-hand sides of these latter two equations are numerically unstable, requiring extended precision arithmetic (at least 45 digits) and about 150 terms in the sums to achieve numerical equality to a reasonable number of significant digits. This suggests the existence of an underlying simplification between the right-hand side sums of which \eqref{a=-I} may be a limiting case - parenthetically, we note that the numerical evaluation of these sums appears to confound both the Maple \cite{Maple} and Mathematica \cite{Math} computer programs.

\section{Summary}

In this work we have summarized the applicability of the general master integration equation \eqref{Master1} to a range of values of a fundamental parameter $a$. Depending on   the value of $a$, extra terms need to be added to \eqref{Master1} as was pointed out in the first Section. In Section 2, these extra terms were located for all possible ranges of the parameter $a$, and in the following Section, explicit examples were presented by employing different values of this parameter. The results were then used to obtain a number of sum rules among the zeros of Riemann's Zeta function in Section 4. This was achieved by using a number of different functions $F(z)$ in \eqref{Master1} together with a judicious choice of the parameter $a$ and application of the principles of the previous sections. It is believed that these sum rules are new, and, represent a template for the derivation of similar rules for a variety of functions using the methods presented here. 	

\appendix
\section{Appendix - Symbols}

Various symbols used in the text are reproduced here (see \cite {Milgram}). Explicit dependence on $\tau=1/2+i\rho$ has been omitted from each, except where necessary to clarify possible ambiguity. From the basic definitions that relate the real and imaginary parts of $\zeta(\tau))$ on the critical line 

\begin{eqnarray}
\zeta_{R}(\tau)=\frac{N}{D_{R}}\zeta_{I}(\tau) \\
\zeta_{I}(\tau)=\frac{N}{D_{I}}\zeta_{R}(\tau)
\label{ZR_ZI}
\end{eqnarray}

we have
\begin{eqnarray}
\nonumber
\mapleinline{inert}{2d}{`D[R]` := 1/2-1/2*(Cp*cos(rhop)+Cm*sin(rhop))/Pi^(1/2);}{%
\[
\mathit{D_{R}} = {\displaystyle \frac {1}{2}}  - {\displaystyle 
\frac {1}{2}} \,{\displaystyle \frac {\mathit{C_{p}}\,\mathrm{cos}(
\mathit{\rho \, \ln(2\pi)}) + \mathit{C_{m}}\,\mathrm{sin}(\mathit{\rho\, \ln(2\pi)})}{\sqrt{
\pi }}} 
\]
}
\\
\nonumber
D_{I}=1-D_{R}
\\
\nonumber
N^2=D_{R}\,D_{I}
\\
\nonumber
\mapleinline{inert}{2d}{Cp := GR*cosh(1/2*Pi*rho)+GI*sinh(1/2*Pi*rho);}{%
\[
\mathit{C_{p}} = \mathit{\Gamma_{R}}\,\mathrm{cosh}({\displaystyle \frac {
\pi \,\rho }{2}} ) + \mathit{\Gamma_{I}}\,\mathrm{sinh}({\displaystyle 
\frac {\pi \,\rho }{2}} )
\]
}
\\
\mapleinline{inert}{2d}{Cm := GI*cosh(1/2*Pi*rho)-GR*sinh(1/2*Pi*rho);}{%
\[
\mathit{C_{m}}= \mathit{\Gamma_{I}}\,\mathrm{cosh}({\displaystyle \frac {
\pi \,\rho }{2}} ) - \mathit{\Gamma_{R}}\,\mathrm{sinh}({\displaystyle 
\frac {\pi \,\rho }{2}} ) 
\]
}
\\
\nonumber
\nonumber
\mapleinline{inert}{2d}{Gr := Re(GAMMA(1/2+i*rho));}{%
\[
\mathit{\Gamma_{R}} = \Re (\Gamma ({\displaystyle \frac {1}{2}}  + i\,
\rho ))
\]
}
\\
\nonumber
\mapleinline{inert}{2d}{Gi := Im(GAMMA(1/2+i*rho));}{%
\[
\mathit{\Gamma_{I}} = \Im (\Gamma ({\displaystyle \frac {1}{2}}  + i\,
\rho ))  \,.
\]
}
\end{eqnarray}
Furthermore,
\begin{eqnarray}
\mapleinline{inert}{2d}{Dp := cos(rho*ln(2*Pi))*GRp+GRm*sin(rho*ln(2*Pi))-Pi^(1/2);}{%
\[
\mathit{D_{+}}= \mathit{\Gamma_{+}}\mathrm{cos}(\rho \,\mathrm{ln}(2\,\pi ))\,
 + \mathit{\Gamma_{-}}\,\mathrm{sin}(\rho \,\mathrm{ln}(2\,
\pi )) - \sqrt{\pi }
\]
}
\\
\nonumber
\mapleinline{inert}{2d}{Dm := GRm*cos(rho*ln(2*Pi))-sin(rho*ln(2*Pi))*GRp+Pi^(1/2);}{%
\[
\mathit{D_{-}} =  - \mathit{\Gamma_{+}}\,\mathrm{sin}(\rho \,\mathrm{ln}(2\,\pi ))\,+\mathit{\Gamma_{-}}\,\mathrm{cos}(\rho \,\mathrm{ln}(2\,
\pi ))
 + \sqrt{\pi }
\]
}
\\
\mapleinline{inert}{2d}{Gp := GR*exp(-1/2*Pi*rho)+exp(1/2*Pi*rho)*GI;}{%
\[
\mathit{\Gamma_{+}} = \mathit{\Gamma_{R}}\,\exp{( - {\pi \,\rho }/{2})} + \mathit{\Gamma_{I}}\,\exp{({\pi \,\rho }/{2})}\,
\]
}
\\
\nonumber
\mapleinline{inert}{2d}{\Gamma_{-} = -GR*exp(1/2*Pi*rho)+exp(-1/2*Pi*rho)*GI;}{%
\[
\mathit{\Gamma_{-}} =  - \mathit{\Gamma_{R}}\,\exp{( {\pi \,\rho }/{2})} + \,\mathit{\Gamma_{I}}\,\exp{
( -  {\pi \,\rho }/{2})}
\]
}
\label{Q5_ids}
\end{eqnarray}
With reference to the alternative notation (see \eqref{Zeta_phase}) introduced in \cite{ReynaDeLune}, it can also be shown after considerable effort, by identifying
\begin{equation}
\frac{\zeta{_R}(\tau)}{\zeta_{I}(\tau)}=-{\displaystyle \cot(\phi(\rho))}
\end{equation}
where \cite[Eq.(9)]{ReynaDeLune} 
\begin{equation}
\phi(\rho)=-\frac{1}{2}\displaystyle{(\gamma+\log\,\pi + 3\,\log\,2+\pi/2)\rho+\sum_{k=0}^{\infty} \left(\frac{2\,\rho}{4k+1}-\arctan(\frac{2\rho}{4k+1}\right)},
\end{equation}
that the two representations are equivalent.
\bibliographystyle{unsrt}
\bibliography{Biblio}


\end{flushleft}
\end{document}